\documentclass[12pt,reqno]{article}

\usepackage{xcolor}
\usepackage{amssymb}
\usepackage{amsmath}
\usepackage{amsthm}
\usepackage[T1]{fontenc}
\usepackage{lmodern}
\usepackage[colorlinks=true,linkcolor=webgreen,filecolor=webbrown,citecolor=webgreen]{hyperref}
\definecolor{webgreen}{rgb}{0,.5,0}
\definecolor{webbrown}{rgb}{.6,0,0}
\usepackage{fullpage}
\usepackage{booktabs,longtable,graphicx,tikz,microtype}

\newcommand{\seqnum}[1]{\href{https://oeis.org/#1}{\rm \underline{#1}}}

\newcommand{\gold}{\alpha}

\theoremstyle{plain}
\newtheorem{theorem}{Theorem}
\newtheorem{corollary}[theorem]{Corollary}
\newtheorem{lemma}[theorem]{Lemma}
\newtheorem{proposition}[theorem]{Proposition}

\theoremstyle{definition}
\newtheorem{definition}[theorem]{Definition}

\theoremstyle{remark}
\newtheorem{remark}[theorem]{Remark}

\hypersetup{pdftitle={Meta-Nested Recurrences},
pdfauthor={Benoit Cloitre, Haobo Ma, Wenlin Zhang}}

\begin{document}

\begin{center}
{\LARGE\bf Meta-Nested Recurrences}
\vskip .3cm
\large
Beno{\^\i}t Cloitre\\
France\\
\href{https://orcid.org/0009-0001-6778-153X}{ORCID: 0009-0001-6778-153X}
\vskip .15cm
Haobo Ma\\
ChronoAI Pte Ltd, Singapore\\
\href{mailto:auric@aelf.io}{\tt auric@aelf.io}
\vskip .15cm
Wenlin Zhang\\
National University of Singapore\\
21 Lower Kent Ridge Road, Singapore 119077, Singapore\\
\href{mailto:e1327962@u.nus.edu}{\tt e1327962@u.nus.edu}
\end{center}

\begin{abstract}
We study recurrences in which a previous sequence value determines the
number of iterations of an inner map. Campbell's example has an explicit
power-of-three formula and inner orbits of eventual period at most two.
Its normalized lower and upper limits are $2/5$ and $3/4$.
For a variable-depth Conway recurrence $C$, we prove the golden-ratio
lower bound, its exact equality set, a Fibonacci upper cap, and the
landing of the inner orbit at every Fibonacci index $F_k$ with $k\ge6$.
These proofs are elementary and start from twenty-four initial values.
The plateau of $C$ to the left of $F_k$ has width exactly
$\lfloor(2k-9)/3\rfloor$, and $C(F_k+t)=F_{k-1}+t$ for
$0\le t\le\lfloor(k-8)/6\rfloor$. A maximum principle turns a finite
exact computation into the bound $\limsup C(n)/n<0.65446923$. We
characterize the convergence of $C(n)/n$ to $\alpha=(\sqrt5-1)/2$ in
terms of negative selected orbit jumps and of a limiting convex profile,
and we derive entropy and variance budgets for the recursive splits.
Along inverse rays of the map $x\mapsto C(x)$, the relative errors are
absolutely summable and the growth is asymptotically Fibonacci. For
Grytczuk's fixed-depth recurrence, we give a Pascal-word proof of
convergence with the sharp error order $n/\sqrt{\log n}$ along interval
maxima. The full variable-depth
limit remains open.
\end{abstract}

\section{Introduction}

Nested recurrences are recurrences in which values at earlier indices occur
inside the arguments of the sequence.  The classical Hofstadter $Q$-sequence
is defined by
\begin{equation}
 Q(n)=Q(n-Q(n-1))+Q(n-Q(n-2)),\qquad Q(1)=Q(2)=1,
 \label{eq:Q}
\end{equation}
and was introduced in \emph{G\"odel, Escher, Bach} \cite{Hofstadter1979}.
The two recurrences below add a second layer, in which the number of
applications of an inner map is itself a previous sequence value.  We call this a
\emph{meta-nested recurrence}.

The surrounding literature contains several Conway-type and Hofstadter-type
recurrences, including the variants studied by Mallows \cite{Mallows1991},
Kubo and Vakil \cite{Kubo1996}, Tanny \cite{Tanny1992}, and Grytczuk
\cite{Grytczuk2004}.  Isgur and Rahman \cite{IsgurRahman2011} study a
fixed-depth family containing Grytczuk's example, prove generalized
Fibonacci anchor identities, and ask whether the normalized sequences
converge.  Section~\ref{sec:grytczuk} proves convergence for their depth-two
case and gives its sharp error order.  The slow-growth property and the
Fibonacci anchors are known. Our proof includes an independent word
derivation, so that the asymptotic argument is self-contained.  Pascal-word
methods and their connection with Conway's sequence also appear in
Janvresse, de la Rue and Velenik \cite{Janvresse2005}.

Fixed iteration counts occur in other families, including Cloitre's
arithmetical solutions \cite{Cloitre2025}.  In
\eqref{eq:campbell} and \eqref{eq:cloitre}, the iteration depth instead
depends on the preceding sequence value.  Cloitre's scaling approach
\cite{CloitreScaling} relates geometric anchor sequences to linear limits
for slow nested recurrences.  For the variable-depth sequence $C$, we prove
an equivalent convergence criterion in terms of selected inner-orbit steps.

The paper has two purposes.  First, it records complete results for the
Campbell example and the global golden structure of the Cloitre example.
Second, it isolates the remaining asymptotic problem in a form that can be
checked locally, namely a negative landing-jump condition in the
interiors of the Fibonacci intervals.  A separate dispersion criterion
gives a conditional rate of convergence.  The proofs of the golden
structure and of the collars are elementary. A finite exact computation
gives a numerical tail bound, and the maximum principle extends it to
all later indices. Along inverse rays of the map $x\mapsto C(x)$, the
normalized sequence converges unconditionally, while the full
variable-depth limit remains open.

\section{Definitions}

\begin{definition}[Campbell recurrence]
Let $b(1)=1$.  For $n\ge2$ put
\begin{equation}
 T_n(x)=n-b(x),\qquad b(n)=T_n^{\,b(n-1)}(n-1),
 \label{eq:campbell}
\end{equation}
where the exponent denotes exactly that many applications of $T_n$.
\end{definition}

\begin{definition}[Cloitre recurrence]
Let $C(1)=C(2)=1$.  For $n\ge3$ put
\begin{equation}
 T_n(x)=n-C(x),\qquad g=T_n^{\,C(n-1)}(n-1),\qquad
 C(n)=C(g)+C(n-g).
 \label{eq:cloitre}
\end{equation}
The prescribed depth is part of the definition. The split is not an
arbitrary periodic point of the inner map.
\end{definition}

We use $F_0=0,F_1=1,F_{j+1}=F_j+F_{j-1}$ and
\[
 \gold=\frac{\sqrt5-1}{2}=\frac1\phi,\qquad
 G(n)=\lfloor\gold(n+1)\rfloor.
\]

\begin{lemma}[Well-definedness and cyclic selection]\label{lem:foundations}
The sequence $C$ is well-defined and satisfies
\[
 \lceil n/2\rceil\le C(n)\le n-1\quad(n\ge2),\qquad
 C(n)\le3n/4\quad(n\ge2),\qquad C(n)\ge3n/5\quad(n\ge3).
\]
For every $n\ge3$, $g(n)$ lies on the cycle reached from $n-1$.
\end{lemma}
The proof is in Appendix~\ref{app:gold-base}. In particular, cycle
reduction computes the prescribed iterate exactly, without changing the
starting point or the iteration depth.

\subsection{Worked examples}\label{subsec:worked-examples}

The definition of $C$ is short, but each value is obtained through a small dynamical computation. We describe this computation for a few values of $n$, each chosen to show one new feature.

Recall that $C(1)=C(2)=1$ and that, for $n\ge 3$,
\[
T_n(y)=n-C(y),\quad h_n=C(n-1),\quad g(n)=T_n^{h_n}(n-1),
\quad C(n)=C(g(n))+C(n-g(n)).
\]
Since $1\le C(y)\le y$, the map $T_n$ sends $\{1,\dots,n-1\}$ into itself and only uses values of $C$ that are already known. We write $y_r=T_n^{r}(n-1)$ for the orbit of $n-1$. This orbit is eventually periodic. Its entry rank $\tau_n$ is the first $r$ for which $y_r$ lies on a cycle, and $p_n$ denotes the length of that cycle. For every $n\ge 3$ one has $h_n\ge\tau_n$ (Lemma~\ref{lem:foundations}). The point $g(n)$ therefore lies on the cycle and
\[
g(n)=y_{\tau_n+((h_n-\tau_n)\bmod p_n)}.
\]
In practice one never iterates $h_n$ times. Once the orbit has reached its cycle, only the residue of $h_n-\tau_n$ modulo $p_n$ matters.

If the depth $h_n$ is replaced by $1$, then $g(n)=n-C(n-1)$ and the recurrence becomes $C(n)=C(n-C(n-1))+C(C(n-1))$. This is the Hofstadter--Conway recurrence $a(n)=a(a(n-1))+a(n-a(n-1))$ of A004001 with its two terms exchanged. The sequence $C$ follows the same map for $C(n-1)$ steps instead of one. Table~\ref{tab:first-values} gives the first values of both sequences together with the data of the computation.

\begin{table}[ht]
\centering
\small
\setlength{\tabcolsep}{3.2pt}
\resizebox{\textwidth}{!}{%
\begin{tabular}{l*{19}{c}}
\toprule
$n$ & 3 & 4 & 5 & 6 & 7 & 8 & 9 & 10 & 11 & 12 & 13 & 14 & 15 & 16 & 17 & 18 & 19 & 20 & 21\\
\midrule
$h_n$ & 1 & 2 & 3 & 3 & 4 & 5 & 5 & 6 & 7 & 7 & 8 & 8 & 9 & 10 & 11 & 12 & 12 & 13 & 13\\
$\tau_n$ & 0 & 0 & 0 & 1 & 3 & 4 & 2 & 2 & 3 & 5 & 6 & 4 & 3 & 5 & 5 & 6 & 4 & 5 & 6\\
$p_n$ & 1 & 2 & 2 & 2 & 1 & 1 & 2 & 2 & 2 & 1 & 1 & 2 & 2 & 2 & 1 & 1 & 2 & 1 & 1\\
$g(n)$ & 2 & 3 & 2 & 3 & 4 & 5 & 5 & 7 & 6 & 7 & 8 & 9 & 8 & 10 & 10 & 11 & 12 & 12 & 13\\
$C(n)$ & 2 & 3 & 3 & 4 & 5 & 5 & 6 & 7 & 7 & 8 & 8 & 9 & 10 & 11 & 12 & 12 & 13 & 13 & 13\\
$a(n)$ & 2 & 2 & 3 & 4 & 4 & 4 & 5 & 6 & 7 & 7 & 8 & 8 & 8 & 8 & 9 & 10 & 11 & 12 & 12\\
\bottomrule
\end{tabular}}
\caption{Depth $h_n=C(n-1)$, entry rank $\tau_n$, period $p_n$, selected point $g(n)$ and value $C(n)$ for $3\le n\le 21$, compared with the Hofstadter--Conway sequence $a$ (A004001). Both sequences start with $1,1$.}
\label{tab:first-values}
\end{table}

\paragraph{The value $n=4$.}
The first value where $C$ and $a$ differ is $n=4$. Since $C(3)=2$ and $C(2)=1$, the orbit of $3$ under $T_4(y)=4-C(y)$ is $3,2,3,2,\dots$, a cycle of length $2$ from the start. Thus $\tau_4=0$, $p_4=2$, and the even depth $h_4=C(3)=2$ brings the orbit back to $y_2=3$. Hence $g(4)=3$ and $C(4)=C(3)+C(1)=3$. A single step stops at $y_1=2$ and gives $C(2)+C(2)=2$, which is the value $a(4)=2$ of the Hofstadter--Conway sequence.

\paragraph{The value $n=10$.}
Here $h_{10}=C(9)=6$. Since $C(9)=6$, $C(4)=3$, $C(7)=5$ and $C(5)=3$, the orbit of $9$ under $T_{10}$ is
\[
9\mapsto 4\mapsto 7\mapsto 5\mapsto 7\mapsto 5\mapsto\cdots.
\]
After a transient of two steps it enters the $2$-cycle $\{5,7\}$, so $\tau_{10}=p_{10}=2$. As $h_{10}-\tau_{10}=4$ is even, the iteration stops at $y_2=7$ and
\[
C(10)=C(7)+C(3)=5+2=7.
\]
An odd depth would have stopped at $5$ and given $C(5)+C(5)=6$. When the cycle has two points, the parity of the depth decides which one is used.

\paragraph{Fibonacci arguments, $n=13$.}
For $n=13$ the depth is $h_{13}=C(12)=8$ and the orbit of $12$ is
\[
12\mapsto 5\mapsto 10\mapsto 6\mapsto 9\mapsto 7\mapsto 8\mapsto 8\mapsto\cdots.
\]
It alternates on both sides of $8$ and reaches the fixed point $8$ at rank $6\le h_{13}$. Therefore $g(13)=8$ and $C(13)=C(8)+C(5)=5+3=8$. The same picture holds at every Fibonacci number $F_k\ge 8$. The value $C(F_{k-1})=F_{k-2}$ makes $F_{k-1}$ a fixed point of $T_{F_k}$, the orbit of $F_k-1$ ends there, and
\[
C(F_k)=C(F_{k-1})+C(F_{k-2})=F_{k-2}+F_{k-3}=F_{k-1}.
\]
The landing statement of Theorem~\ref{thm:cap} shows that this happens for every $k\ge 6$. The value $n=5$ behaves differently, since the orbit of $4$ is the $2$-cycle $\{2,4\}$ and the odd depth $h_5=3$ selects $2$, which still gives $C(5)=C(2)+C(3)=3$.

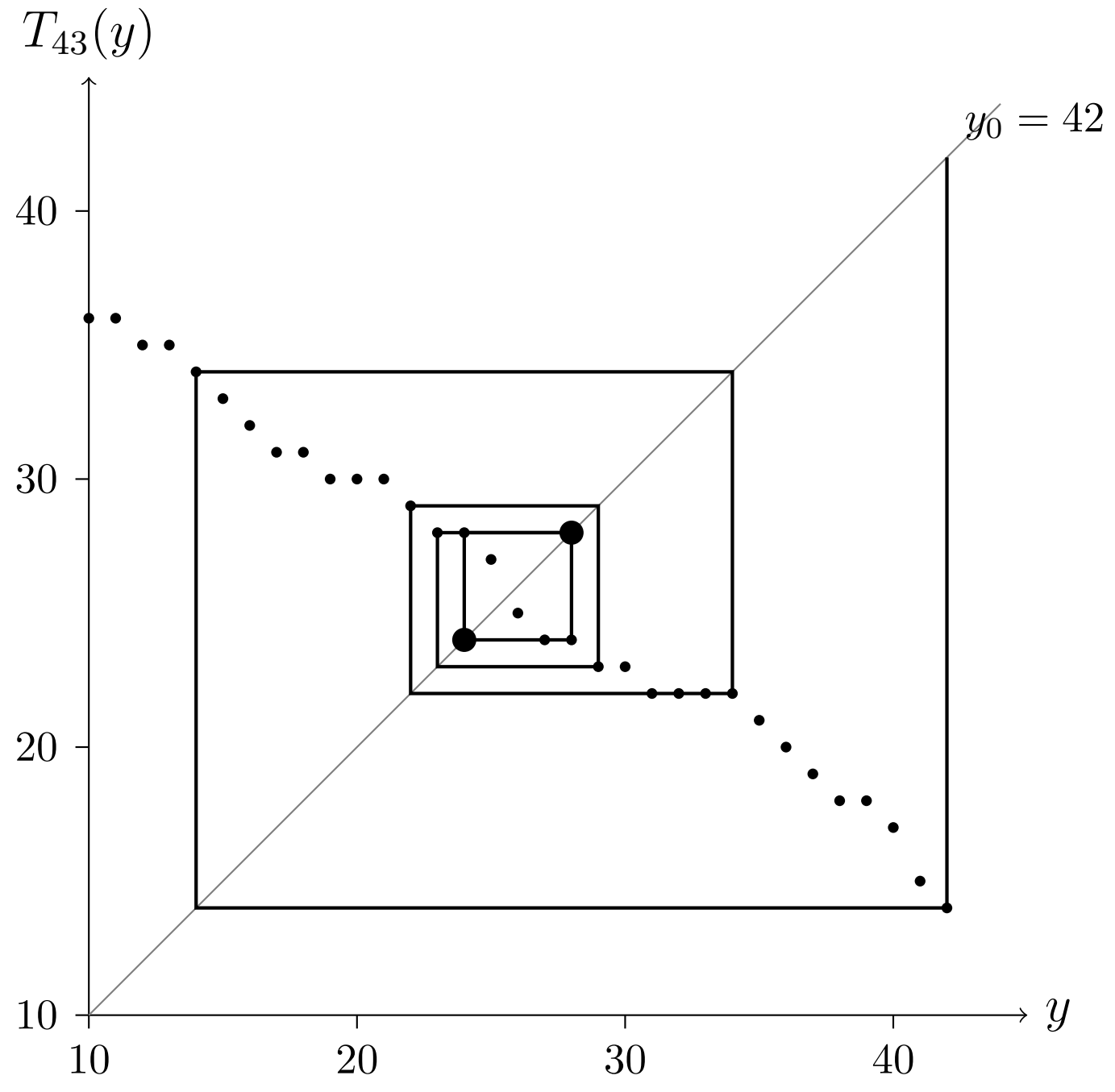
\begin{figure}[ht]
\centering
\begin{tikzpicture}[x=0.22cm,y=0.22cm]
  \draw[->] (10,10) -- (45,10) node[right] {$y$};
  \draw[->] (10,10) -- (10,45) node[above] {$T_{43}(y)$};
  \foreach \x in {10,20,30,40} {
    \draw (\x,10) -- (\x,9.5) node[below] {\footnotesize $\x$};
    \draw (10,\x) -- (9.5,\x) node[left] {\footnotesize $\x$};
  }
  \draw[gray] (10,10) -- (44,44);
  \foreach \p in {(10,36),(11,36),(12,35),(13,35),(14,34),(15,33),(16,32),(17,31),(18,31),(19,30),(20,30),
                  (21,30),(22,29),(23,28),(24,28),(25,27),(26,25),(27,24),(28,24),(29,23),(30,23),(31,22),
                  (32,22),(33,22),(34,22),(35,21),(36,20),(37,19),(38,18),(39,18),(40,17),(41,15),(42,14)}
    \fill \p circle (0.2);
  \draw[thick] (42,42) -- (42,14) -- (14,14) -- (14,34) -- (34,34) -- (34,22) -- (22,22) -- (22,29)
               -- (29,29) -- (29,23) -- (23,23) -- (23,28) -- (28,28) -- (28,24) -- (24,24) -- (24,28) -- (28,28);
  \node[above right] at (42,42) {\footnotesize $y_0=42$};
  \fill (24,24) circle (0.45) (28,28) circle (0.45);
\end{tikzpicture}
\caption{The orbit of $42$ under $T_{43}(y)=43-C(y)$. The small dots are the points $(y,T_{43}(y))$ for $10\le y\le 42$, the grey line is the diagonal, and the two large dots mark the $2$-cycle $\{24,28\}$.}
\label{fig:cobweb43}
\end{figure}

\paragraph{The first descent, $n=43$.}
The sequence $C$ is not monotone, and its first descent shows the role of the depth. For $n=43$ the depth is $h_{43}=C(42)=29$. The orbit of $42$ is listed in Table~\ref{tab:orbit43} and drawn in Figure~\ref{fig:cobweb43}. It alternates on both sides of the diagonal with gaps $28,20,12,7,6,5,4$ and enters the $2$-cycle $\{24,28\}$ at rank $\tau_{43}=6$. Since $h_{43}-\tau_{43}=23$ is odd, the iteration stops at $y_7=24$, and
\[
C(43)=C(24)+C(19)=15+13=28<29=C(42).
\]
The other point of the cycle would give $C(28)+C(15)=19+10=29$, and a single step would also give $29$, since $y_1=14$ and $C(14)+C(29)=9+20=29$. The descent at $n=43$ comes from the parity of the depth alone. The first increment larger than $1$ occurs earlier, at $n=26$, where $C(25)=16$ and $C(26)=18$.

\begin{table}[ht]
\centering
\begin{tabular}{l*{8}{c}}
\toprule
$r$ & 0 & 1 & 2 & 3 & 4 & 5 & 6 & 7\\
\midrule
$y_r$ & 42 & 14 & 34 & 22 & 29 & 23 & 28 & 24\\
$C(y_r)$ & 29 & 9 & 21 & 14 & 20 & 15 & 19 & 15\\
$y_{r+1}=43-C(y_r)$ & 14 & 34 & 22 & 29 & 23 & 28 & 24 & 28\\
\bottomrule
\end{tabular}
\caption{The orbit of $42$ under $T_{43}$. From rank $6$ on it alternates between $28$ and $24$.}
\label{tab:orbit43}
\end{table}

\paragraph{Longer cycles.}
Cycles of length $1$ and $2$ are the only ones for $n<68$. The first cycles of length $4$, $6$, $3$ and $5$ appear at $n=68$, $120$, $126$ and $196$. For $n=196$ the orbit of $195$ enters the $5$-cycle $(120,116,117,118,115)$ at rank $8$. The depth is $h_{196}=C(195)=131$ and $131-8\equiv 3\pmod 5$, so the iteration stops three points after the entry, at $y_{11}=118$, and $C(196)=C(118)+C(78)=81+53=134$. The five points of this cycle would give the five values $129,132,131,134,130$.

\section{Campbell's recurrence}

\begin{theorem}[Explicit formula]\label{thm:campbell}
For every $n\ge2$, let $s$ be the unique power of three in the indicated
range.  Then
\begin{equation}
 b(n)=
 \begin{cases}
 \min(n-s,3s),& n\text{ even and }2s\le n<6s,\\
 \max(2s,n-3s),& n\text{ odd and }3s\le n<9s.
 \end{cases}
 \label{eq:campbell-formula}
\end{equation}
\end{theorem}

\begin{proof}
The twelve orbit rows in Appendix~\ref{app:campbell} give a strong
induction proof, including the parity of the prescribed depth. In
particular, the orbit agrees with its two-periodic tail from rank $4$.
Even depth selects rank $4$, and odd depth selects rank $5$.
\end{proof}

\begin{corollary}
The recurrence \eqref{eq:campbell} is globally well-defined.  Its inner
orbit has transient length at most four before entering a fixed point or a
two-cycle, and
\[
 b(3n)=3b(n)\quad(n\ge2),\qquad
 \liminf_{n\to\infty}\frac{b(n)}n=\frac25,
 \qquad
 \limsup_{n\to\infty}\frac{b(n)}n=\frac34.
\]
In particular, $b(n)/n$ does not converge.
\end{corollary}

\section{The global golden structure of \texorpdfstring{$C$}{C}}

The main comparison sequence is the Hofstadter $G$-sequence $G(n)$.  Its
identity $G(n)=n-G(G(n-1))$ follows directly from the rotation identity for
$\gold$ and uniquely determines the sequence from its initial values.

\begin{theorem}[Golden lower bound and equality set]\label{thm:golden}
For every $n\ge1$,
\[
 C(n)\ge G(n).
\]
Moreover, equality holds exactly on
\begin{equation}
 Z=\{F_j,F_j+1:j\ge2\}
   \cup\{F_j-1:j\ge3\text{ odd}\}
   \cup\{11,24,25,59\}.
 \label{eq:equality-set}
\end{equation}
\end{theorem}

\begin{theorem}[Fibonacci upper cap and landing]\label{thm:cap}
If $F_j\le n<F_{j+1}$ and $j\ge3$, then
\begin{equation}
 C(n)\le U(n):=\min(n-F_{j-2},F_j).
 \label{eq:upper-cap}
\end{equation}
Consequently,
\begin{equation}
 \begin{aligned}
 C(F_j)&=F_{j-1}&& (j\ge2),\\
 C(F_j-1)&=F_{j-1}&& (j\ge5),\\
 C(F_j+1)&=F_{j-1}+1&& (j\ge3).
 \end{aligned}
 \label{eq:fib-values}
\end{equation}
For $k\ge6$, every orbit of $T_{F_k}$ on $[1,F_k-1]$ eventually reaches
the fixed point $F_{k-1}$.  In particular, the prescribed depth in
\eqref{eq:cloitre} selects that point.
\end{theorem}

The cap $U$ has the same block formula as Peter Bala's sequence
\seqnum{A356988}. The two block descriptions in its OEIS entry
\cite{OEIS} meet at the Lucas number $F_j+F_{j-2}$.
This identification rests on that block formula, not on the nested
recurrence that defines A356988.

\begin{corollary}
\begin{equation}
 \liminf_{n\to\infty}\frac{C(n)}n=\gold.
 \label{eq:liminf}
\end{equation}
The full limit is not established by these results.
\end{corollary}

\begin{proof}[Proof of Theorems~\ref{thm:golden} and~\ref{thm:cap} and of the corollary]
Appendix~\ref{app:golden} gives a simultaneous induction for the lower
bound, for the implication $C(x)=G(x)\Rightarrow x\in Z$ and for the
cap, strengthened at the changes of slope of $U$. It displays the twenty-four initial values and all
exceptional child pairs. The reverse inclusion in the equality set
and the Fibonacci landing statement then follow from the comparison
and capture bounds. The corollary follows from the lower bound and the
anchors $C(F_j)=F_{j-1}$.
\end{proof}

\begin{theorem}[The left Fibonacci plateau]
\label{article:left-plateau}
For $k\ge8$, put $w_k=\lfloor(2k-9)/3\rfloor$. Then
\[
 C(F_k-u)=F_{k-1}
 \quad\Longleftrightarrow\quad
 0\le u\le w_k
 \qquad(0\le u\le F_{k-2},\ u\in\mathbb Z).
\]
Moreover,
\[
 C(F_k-w_k-1)=F_{k-1}-1,
 \qquad
 \min\{n\ge1\mid C(n)=F_{k-1}\}=F_k-w_k.
\]
\end{theorem}

The width $w_k$ increases by one exactly when $F_{k-1}$ is odd, so it
grows linearly in $k$. The proof is an induction on $k$ and is given in
Appendix~\ref{app:plateau-collar}.

\section{Fibonacci collars and profile renormalization}

For any $j\ge2$ and integer $u$ with $F_j+u\ge1$, define
\[
 P_j(u)=C(F_j+u)-F_{j-1},\qquad
 D_j(u)=P_j(u)-\max(0,u),\qquad E_j(u)=-D_j(u).
\]
For the capture theorem, fix $k\ge6$, write $n=F_k+t\ge3$ with
$t\in\mathbb Z$, and put $A=F_{k-1}$.

\begin{theorem}[Collar capture and profile identity]\label{thm:capture}
Every cycle of $T_n$ is contained in
\begin{equation}
 I_{k,t}=[A+\min(0,t),A+\max(0,t)],
 \label{eq:collar}
\end{equation}
and every orbit eventually enters this interval.  Hence every cycle has
period at most $|t|+1$.  If the selected split is $a=A+s$, then
\begin{equation}
 \begin{aligned}
 P_k(t)&=P_{k-1}(s)+P_{k-2}(t-s),\\
 D_k(t)&=D_{k-1}(s)+D_{k-2}(t-s)\le0,
 \end{aligned}
 \label{eq:profile}
\end{equation}
with $s\in[\min(0,t),\max(0,t)]\cap\mathbb Z$.
\end{theorem}

\begin{proof}
Put $B=F_{k-2}$, so $n-B=A+t$, and write $I=[\ell,h]$ for
\eqref{eq:collar}. The cones~\eqref{eq:strict-cones} hold at the
anchor $(A,B)$ for every actual argument $1\le x<n$. They make
$I\cap[1,n-1]$ invariant. For $u\ge1$, they also give
\[
 T_n(\ell-u)\in[\ell,h+u-1],\qquad
 T_n(h+u)\in[\ell-u,h]
\]
whenever the argument belongs to the domain. An exterior point therefore
enters $I$ or switches sides, with integer distance decreasing within
at most two steps. This proves capture for either sign of $t$, even if
the anchor lies outside the domain. A cycle has distinct integer points,
so its period is at most $|t|+1$.

The selected split is cyclic by Lemma~\ref{lem:foundations}, hence
$a=A+s$ with $\min(0,t)\le s\le\max(0,t)$. Its complement is
$n-a=F_{k-2}+(t-s)$. Both are positive actual arguments.
Substitution in \eqref{eq:cloitre} gives the first profile identity.
The offsets $s$ and $t-s$ have the same weak sign as $t$, so
\[
 \max(0,s)+\max(0,t-s)=\max(0,t),
\]
which gives the identity for $D_k$.
The global cap $U$ of Theorem~\ref{thm:cap}, with $U(1)=1$, is nondecreasing,
has increments zero or one, and satisfies $U(F_j)=F_{j-1}$.
Consequently, for every integer $u$ with $F_j+u\ge1$,
\[
 C(F_j+u)\le U(F_j+u)\le F_{j-1}+\max(0,u).
\]
Applying this bound to the two child offsets gives their nonpositive
defects and completes the proof.
\end{proof}

\begin{corollary}[Sublinear-width convergence]
For arbitrary integer offsets with $t_k/F_k\to0$, including negative offsets,
\[
 \frac{C(F_k+t_k)}{F_k+t_k}\longrightarrow\gold.
\]
For every fixed integer $t$, the periods are uniformly bounded in $k$.
For every fixed integer $t$, for all sufficiently large $k$,
\[
 C(F_k+t)=F_{k-1}+\max(0,t).
\]
\end{corollary}

For the first statement the lower bound and the cap suffice, since on
either side of $F_k$ the cap differs from $F_{k-1}$ by at most
$|t_k|+O(1)$. The left plateau of Theorem~\ref{article:left-plateau}
and the right collar of Proposition~\ref{article:log-collar} give the
fixed-offset identity.

For each fixed integer $t$, applying this identity at level $k-1$ to the
finitely many offsets $\min(0,t)\le s\le\max(0,t)$ gives, for all
sufficiently large $k$,
\[
 T_{F_k+t}(F_{k-1}+s)=F_{k-1}+t-\max(0,s).
\]
Thus, if $t<0$, every orbit reaches the unique fixed point $F_{k-1}+t$;
if $t\ge0$, the captured map is the reflection $s\mapsto t-s$, so every
cycle has period one or two.

For the five-cycle at $196$ in Section~\ref{subsec:worked-examples},
the offsets from $F_{11}=89$ are $31\to27\to28\to29\to26\to31$.
The profile values $E_{11}$ along it are $(6,3,5,3,5)$.

\begin{proposition}[A logarithmic right collar]
\label{article:log-collar}
For every $k\ge8$,
\[
 C(F_k+t)=F_{k-1}+t
 \qquad
 \left(0\le t\le\left\lfloor\frac{k-8}{6}\right\rfloor,
 \ t\in\mathbb Z\right).
\]
\end{proposition}

The width grows logarithmically with $F_k$. The proof in
Appendix~\ref{app:plateau-collar} shows that a first positive defect
would propagate to the levels $k-2$ and $k-4$ and force
$F_{k-1}$, $F_{k-3}$ and $F_{k-5}$ to have the same parity, whereas
exactly one of them is even.

For $0\le t\le F_{k-1}$, write a captured point as $F_{k-1}+u$. The inner map in these coordinates is
\begin{equation}
 u\longmapsto t-u+E_{k-1}(u).
 \label{eq:reflection}
\end{equation}
For a five-cycle this gives
\begin{equation}
 2u_0=t+E_0-E_1+E_2-E_3+E_4,
 \qquad
 2\sum_{i=0}^4u_i=5t+\sum_{i=0}^4E_i.
 \label{eq:fiveclosure}
\end{equation}
Here $E_i=E_{k-1}(u_i)$. These closure identities hold for every
five-cycle. They do not classify the cycles, since the defect word and
the phase selected by the prescribed start and depth remain to be
determined.

\section{An interior landing-jump criterion}

Write $g(n)=T_n^{C(n-1)}(n-1)$ and define
\begin{equation}
 \begin{aligned}
 d(n)&=\frac{C(n)}n-\alpha,\qquad e_-(n)=\max(0,-e(n)),\\
 e(n)&=T_n(g(n))-g(n)=n-C(g(n))-g(n).
 \end{aligned}
 \label{eq:landing-error}
\end{equation}
The displacement $e(n)$ is an inner-orbit step, whereas $d(n)$ measures
the outer sequence's relative error.  For $F_j\le x<F_{j+1}$ put
\[
 q(x)=\frac{x-F_j}{F_{j-1}}\in[0,1).
\]
For $0<\delta<1/2$, define the nonnegative localized limsup
\[
 E_-(\delta)=\limsup_{n\to\infty}
 \mathbf1_{\{\delta\le q(g(n))\le1-\delta\}}\frac{e_-(n)}n.
\]
This definition gives zero if only finitely many indices meet the restriction.

Cloitre \cite[Theorem~3]{CloitreScaling} derives linear limits from geometric
anchors for slow recurrences under bounded-error splitting hypotheses.
His Proposition~24 treats an alternately
perturbed Conway recurrence under geometric-anchor and sublinear
adjacent-jump hypotheses.
Neither set of hypotheses has been established for the variable-depth
sequence $C$ here.  Our criterion instead concerns only the negative part
of the prescribed inner landing step, and it does not require a bound on
every adjacent difference $C(n)-C(n-1)$.

By Lemma~\ref{lem:foundations}, the selected point is cyclic. The
capture bounds~\eqref{eq:child-rectangle} give, with inherited endpoint labels,
\begin{equation}
 F_{k-1}\le g(n)\le F_k,\qquad
 F_{k-2}\le n-g(n)\le F_{k-1}
 \quad(F_k\le n<F_{k+1},\ n\ge8).
 \label{eq:closed-children}
\end{equation}

\begin{theorem}\label{thm:landing-criterion}
Put $L=\limsup_{n\to\infty}d(n)$ and
$E_-=\limsup_{n\to\infty}e_-(n)/n$.  Then
\begin{equation}
 0\le L\le\phi E_-.
 \label{eq:one-sided}
\end{equation}
The following conditions are equivalent.
\begin{enumerate}
\item $C(n)/n\longrightarrow\alpha$.
\item $e_-(n)=o(n)$.
\item $E_-(\delta)=0$ for every $\delta\in(0,1/2)$.
\end{enumerate}
If $L>0$, there is a $\delta>0$ such that $E_-(\delta)\ge\alpha L$.
\end{theorem}

\begin{proof}
The exact splitting identity gives
\begin{equation}
 d(n)=\frac{g(n)}n d(g(n))+
        \frac{n-g(n)}n d(n-g(n)).
 \label{eq:extremal-average}
\end{equation}
Both weights stay bounded away from zero by \eqref{eq:closed-children}.
Hence, along any subsequence with $d(n)\to L$, both child errors tend to
$L$, since each is at most $L+o(1)$ and a fixed deficit in either child
would force a fixed deficit in the parent.

Suppose $L>0$.  Set $\theta(n)=n/F_k$ when $F_k\le n<F_{k+1}$.
Among the cluster points $(L,\theta)$ of $(d(n),\theta(n))$, choose one
with minimum second coordinate $\theta_*$.  Take a subsequence realizing
it, and pass further so that $g(n)/n\to\lambda$.  Since $d(g(n))\to L>0$
but Fibonacci endpoints have error tending to zero, the child lies strictly
between $F_{k-1}$ and $F_k$ eventually.  Its canonical phase is therefore
$g(n)/F_{k-1}$, and minimality yields
\[
 \theta_*\le\lim\frac{g(n)}{F_{k-1}}
       =\lambda\phi\theta_*,\qquad \lambda\ge\alpha.
\]
Using \eqref{eq:landing-error}, we obtain
\[
 \lim\frac{e(n)}n=1-\lambda(1+\alpha+L)
       \le1-\alpha(\phi+L)=-\alpha L.
\]
This proves \eqref{eq:one-sided}.

The cap gives
$C(x)\le\min(F_{j-1}+x-F_j,F_j)$ on $[F_j,F_{j+1})$.
Writing $x=F_j+qF_{j-1}$ and using the Fibonacci ratios, it gives,
uniformly in $q\in[0,1]$,
\begin{equation}
 d(x)\le H(q)+o(1),\qquad
 H(q)=\frac{\min(\alpha^3q,\alpha^2(1-q))}{1+\alpha q}.
 \label{eq:phase-cap}
\end{equation}
Since $H$ is continuous and vanishes at both endpoints, $d(g(n))\to L>0$
forces $q(g(n))\in[\delta,1-\delta]$ eventually for some fixed $\delta>0$.
The same subsequence thus proves $E_-(\delta)\ge\alpha L$.
Condition (3) consequently forces $L=0$, and the golden lower bound gives (1).

To prove (1) implies (2), let $u,v$ be the minimum and maximum over all
cycle points of $T_n$.  Each has a cycle predecessor.  Cycle capture
keeps their arguments at a positive fraction of $n$.  If (1) holds, then,
for any small $\varepsilon>0$ and all large $n$, all these arguments obey
$mx\le C(x)\le rx$, where $m=\alpha-\varepsilon$ and
$r=\alpha+\varepsilon<1$.  Taking predecessors gives
\[
 u\ge n-rv,\qquad v\le n-mu,
 \qquad v-u\le\frac{r-m}{1-mr}n.
\]
Letting $\varepsilon\downarrow0$ shows that every cycle has width $o(n)$,
including every selected step.  Therefore $e_-(n)=o(n)$.
Finally, (2) immediately implies (3).
\end{proof}

The phase minimum selects an extremal child with a macroscopic negative
jump. The cap then keeps that child away from both Fibonacci endpoints.
Thus it remains to estimate negative jumps in each fixed interior interval.

\begin{corollary}
For any fixed $0<\delta\le\alpha^2$,
\[
 L\le\max\left\{\phi E_-(\delta),
        \frac{\alpha^3\delta}{1-\alpha^2\delta}\right\}.
\]
\end{corollary}
\begin{proof}
On the two boundary strips $[0,\delta]\cup[1-\delta,1]$, the maximum of
$H$ in \eqref{eq:phase-cap} is
$\alpha^3\delta/(1-\alpha^2\delta)$.  If $L$ exceeds this number,
the extremal child subsequence must lie in the chosen interior, giving
$L\le\phi E_-(\delta)$.  Otherwise the second bound applies.
\end{proof}

\begin{remark}
If $p$ is the actual cycle predecessor of $g$, then
$n-g=C(p)$ and $e(n)=C(p)-C(g)$.  In a two-cycle,
$C(g)-C(p)=g-p$. Equivalently, with $D(x)=C(x)-\alpha x$,
$D(g)-D(p)=\alpha^2(g-p)$.  This turns a negative jump on that two-cycle
into a secant of slope one.  Longer cycles require their own argument,
since the predecessor need not be $T_n(g)$.
\end{remark}

Exact computation through $2^{28}$ gives the maximum of $e_-(n)/n$ on
each dyadic shell $[2^m,2^{m+1})$. On $[2^{27},2^{28})$ it is
$653668/166238801\approx0.003932102$, attained at $n=166238801$ with
$g=102990319$ and $T_n(g)=102336651$, and on $[2^{23},2^{24})$ it is
$45807/9274448\approx0.004939054$. The shell maxima are not monotone
before $2^{21}$ and decrease slowly from there on. From $2^{21}$ on, each
of them is attained on a two-cycle just to the right of a Fibonacci
number. These finite data support investigation of the criterion but
do not bound its limiting limsup.

\section{A maximum principle and entropy budgets}
\label{article:maximum-entropy}

We retain the canonical split $g(n)=T_n^{C(n-1)}(n-1)$ and write
\[
 d(n)=C(n)/n-\alpha,\qquad
 e(n)=n-g(n)-C(g(n)),\qquad L=\limsup d(n).
\]
All child labels below are inherited from closed Fibonacci blocks.
The comparison and capture results established earlier are used
throughout this section.

\begin{proposition}[Maximum principle]
For $j\ge8$,
\begin{equation}
\label{article:maximum}
 \sup_{n\ge F_j}\frac{C(n)}n
 =\max_{F_j\le n<F_{j+2}}\frac{C(n)}n.
\end{equation}
\end{proposition}
\begin{proof}
For $n\ge F_{j+2}\ge F_{10}=55$, Lemma~\ref{lem:foundations}
and capture put both selected children in $[F_j,n)$.
Their ratios average to the parent ratio with weights $g(n)/n$ and
$(n-g(n))/n$. Strong induction therefore extends the finite maximum
to the whole tail.
\end{proof}

An exact computation through $F_{45}$ gives the finite certificate
\[
 561724478\,C(n)\le367631381\,n\qquad(F_{43}\le n<F_{45}),
\]
with equality only at $n=561724478$. The accompanying code contains
a hash-table evaluator, a separate array checker and an independent
Brent evaluator. The latter regenerated the full prefix through $F_{45}$
from the initial values and confirmed the inequality and unique equality
on all $701408733$ indices of $[F_{43},F_{45})$ by integer cross products.
The two evaluators give identical maxima on every complete Fibonacci
block from order $20$ through $44$.
Thus the same inequality holds for every $n\ge F_{43}$, and
\[
 \limsup C(n)/n\le367631381/561724478<0.65446923.
\]
The finite verification is a computational premise. Its extension
to an infinite tail is \eqref{article:maximum}.

For $0\le t\le F_{k-1}$, define
\[
 E_k(t)=-D_k(t)=F_{k-1}+t-C(F_k+t),\qquad q=t/F_{k-1}.
\]
As before, $E_k=-D_k$ is the nonnegative profile defect and $e(n)$ is
the selected orbit jump.
Let $f_k$ interpolate the grid values $E_k(t)/F_{k-1}$, and for
$j\ge8$ let $\Psi_j$ be the greatest common convex minorant of
$f_j,f_{j+1}$.
If $F_k+t\ge47$ and $g(F_k+t)=F_{k-1}+s$, Theorem~\ref{thm:capture} and
the child localization~\eqref{eq:child-rectangle} give
\[
 E_k(t)=E_{k-1}(s)+E_{k-2}(t-s),\qquad
 0\le s\le F_{k-2},\quad0\le t-s\le F_{k-3}.
\]

\begin{proposition}[Convex envelopes]
The functions $\Psi_j$ increase and converge uniformly to a continuous convex
function $\Psi_\infty$ satisfying
\[
 (q-\alpha)_+\le\Psi_\infty(q)\le\alpha^2q,
 \qquad \Psi_\infty(0)=0,\quad\Psi_\infty(1)=\alpha^2.
\]
Moreover,
\begin{equation}
\label{article:envelope-formula}
 L=\max_{0\le q\le1}
 \frac{\alpha^2q-\Psi_\infty(q)}{\phi+q}.
\end{equation}
In particular, convergence of $C(n)/n$ to $\alpha$ is equivalent
to $\Psi_\infty(q)=\alpha^2q$ on $[0,1]$.
\end{proposition}
\begin{proof}
Put $a_k=F_{k-2}/F_{k-1}$ and $b_k=1-a_k$.
The child phases $q_0=s/F_{k-2}$ and $q_1=(t-s)/F_{k-3}$ satisfy
\[
 q=a_kq_0+b_kq_1,\qquad
 f_k(q)=a_kf_{k-1}(q_0)+b_kf_{k-2}(q_1).
\]
Jensen's inequality propagates any common convex minorant of two
successive profiles to every later profile. Thus $\Psi_j\le\Psi_{j+1}$.
The cap gives the common lower bound
$\bigl(q-\max\{a_j,a_{j+1}\}\bigr)_+$.
The endpoint values of $\Psi_j$ are zero and
$\min\{F_{j-3}/F_{j-1},F_{j-2}/F_j\}<\alpha^2$, so convexity
gives the upper bound. These bounds imply continuity of the limit
at the endpoints. Convexity gives interior continuity, and monotone
convergence on the compact interval is uniform.

For every real $M$, the affine function
$\alpha^2q-M(\phi+q)$ minorizes both $f_j,f_{j+1}$ exactly
when it minorizes $\Psi_j$. Hence the maximum of
$(\alpha^2q-f_k(q))/(\phi+q)$ over those two profiles equals
the corresponding maximum with $\Psi_j$. A quotient of affine
functions on an interpolation segment attains its maximum at an
endpoint. At the grid points,
\[
 d(F_k+t)=
 \frac{\alpha^2q-f_k(q)+(-1)^k\alpha^k/F_{k-1}}
 {F_k/F_{k-1}+q}.
\]
The discrepancy from the preceding quotient is uniform and tends
to zero. The maximum principle and uniform envelope convergence
prove \eqref{article:envelope-formula}.
\end{proof}

The offset walk chooses the two children with probabilities
$\sigma=s/t$ and $1-\sigma$ when $t>0$. A zero offset has zero
probability. Put $M=E_k(t)/t$ and $\lambda=M/\alpha^2$.
Both are martingales. A root is the starting node of this walk,
with values $q_{\mathrm{init}}$, $M_{\mathrm{init}}$ and
$\lambda_{\mathrm{init}}$. A stopping time is admissible if every
active node before it has size at least $47$.
For a child of zero offset, assign its ratio $M^{(i)}$ any value
in $[0,1]$. With child ratios $M^{(0)},M^{(1)}$, define
\[
 K=\operatorname{KL}((\sigma,1-\sigma)\Vert(a_k,b_k)),\qquad
 V=\sigma(1-\sigma)(M^{(0)}-M^{(1)})^2.
\]

\begin{proposition}[Offset entropy and a residual]
For every bounded admissible stopping time $T$,
\begin{align*}
 \mathbb E\sum_{h<T}K_h
 &=\mathbb E\log q_T-\log q_{\mathrm{init}}
 \le\log(1/q_{\mathrm{init}}),\\
 \mathbb E\sum_{h<T}V_h
 &=\mathbb E M_T^2-M_{\mathrm{init}}^2
 \le M_{\mathrm{init}}(1-M_{\mathrm{init}}).
\end{align*}
The exact jump identity and the residual identity are
\begin{align*}
 e(F_k+t)&=t-2s+E_{k-1}(s),\\
 R&=e/t+2a_k-1-a_kM\\
  &=(M-2)(\sigma-a_k)+\sigma(1-\sigma)(M^{(0)}-M^{(1)}).
\end{align*}
They imply
\[
 R^2\le4K+V/2,\qquad
 \mathbb E\sum_{h<T}R_h^2
 \le4\log(1/q_{\mathrm{init}})
 +M_{\mathrm{init}}(1-M_{\mathrm{init}})/2.
\]
\end{proposition}
\begin{proof}
The child phases are $q\sigma/a_k$ and $q(1-\sigma)/b_k$.
Conditioning gives the logarithmic drift $K$. The conditional
variance of $M$ is $V$. Telescoping and $q\le1$, $M^2\le M$
give the two budgets. The jump formula follows by substituting
$C(g)=F_{k-2}+s-E_{k-1}(s)$ in the definition of $e$.
Subtracting the child average for $M$ gives the residual identity.
Since $|M-2|\le2$ and $\sigma(1-\sigma)\le1/4$,
$R^2\le8(\sigma-a_k)^2+V/2$. The second derivative of the binary
divergence is $1/(z(1-z))\ge4$, giving $K\ge2(\sigma-a_k)^2$.
\end{proof}

These budgets constrain how a profile deficit can persist. Suppose roots at
levels $k_i$ satisfy $q_{\mathrm{init}}\ge a>0$ and
$\lambda_{\mathrm{init}}\le\theta<1$.
Take $H_i\to\infty$ with $k_i-2H_i\to\infty$, fix
$\theta'\in(\theta,1)$, and put $c=\alpha^3(1-\theta')/2$.
For every $b\in(0,1)$,
\begin{equation}
\label{article:occupation}
 \liminf_i\frac1{H_i}\sum_{h<H_i}
 \mathbb P\{q_h\ge b,\ e(x_h)\le-ct_h\}
 \ge1-\frac\theta{\theta'}-
 \frac{\log(1/a)}{\log(1/b)}.
\end{equation}
Indeed, Markov's inequality controls the events $\lambda_h>\theta'$
and $q_h<b$. The residual budget makes $R_h$ tend to zero in
Ces\`aro probability, while uniformly on these trees
\[
 R_h=e(x_h)/t_h+\alpha^3(1-\lambda_h)+o(1).
\]
The claimed negative jump follows on the remaining mass. Both
endpoints of that jump lie on the selected cycle, so its diameter
is at least $ct_h$. This includes every cycle period.

For completeness, on the grid $q=t/F_{k-1}$ with
$1\le t\le F_{k-1}$, let $B_k(\theta)$ be the largest coordinate
with $E_k(t)/(\alpha^2t)\le\theta$, and let $J_k(c)$ be the largest
coordinate with $e(F_k+t)\le-ct$. Both maxima are zero for
an empty set. Theorem~\ref{thm:landing-criterion} and
$|e(F_k+t)|\le t$ give
\[
 C(n)/n\to\alpha
 \quad\Longleftrightarrow\quad
 B_k(\theta)\to0\ \ (0<\theta<1)
 \quad\Longleftrightarrow\quad
 J_k(c)\to0\ \ (c>0).
\]
For the first equivalence, use the exact formula for $d(F_k+t)$
above away from $q=0$ and the cap near $q=0$. For the second,
convergence forces all cycle diameters to be $o(n)$. Conversely,
if $J_k(c)\to0$, then either $e_-(n)<ct$ or
$e_-(n)/n\le t/n\le J_k(c)$.
At $t=0$, the Fibonacci landing statement of Theorem~\ref{thm:cap} gives
$e(F_k)=0$ for every $k\ge6$.
Letting $c$ tend to zero proves $e_-(n)=o(n)$, and
Theorem~\ref{thm:landing-criterion} applies.
The convergence assertion remains a conjecture.

\section{Martingales and the dispersion route}

Use labelled states $(j,N)$ with $F_j\le N\le F_{j+1}$, and let the
selected split be $N=a+b$ with inherited child labels $j-1,j-2$.
Retain these labels at Fibonacci endpoints rather than recomputing a canonical
order.  Choose the first child with probability $a/N$ and the second with
probability $b/N$.  Define
\[
 Z(N)=\frac{C(N)-\gold N}{N},\qquad X(j,N)=\frac{F_j}{N}.
\]

\begin{proposition}[Exact martingales and local variance]
Both $Z$ and $X$ are martingales for this size-biased child choice.  If
$A=F_{j-1}$ and $B=F_{j-2}$, the one-step conditional variance of $X$ is
\begin{equation}
 v(j,N)=\frac{(Ab-Ba)^2}{N^2ab}.
 \label{eq:variance}
\end{equation}
Writing $V_0=0$, the accumulated variance satisfies
\begin{equation}
 V_m(j,N)=v(j,N)+\frac aN V_{m-1}(j-1,a)
                    +\frac bN V_{m-1}(j-2,b),
 \label{eq:variance-rec}
\end{equation}
with terminal orders stopped at $4$ or $5$.
\end{proposition}
\begin{proof}
The martingale identities follow from $C(N)=C(a)+C(b)$ and
$F_j=F_{j-1}+F_{j-2}$. Subtracting the weighted mean of the two
child values of $X$ gives \eqref{eq:variance}, and conditioning on the
first child choice gives \eqref{eq:variance-rec}.
\end{proof}

\begin{theorem}[Conditional quartic decay]
Assume that there are $\kappa>0$ and $J$ such that every labelled state
$(j,N)$ reached by the selected splits satisfies
\begin{equation}
 V_4(j,N)\ge \kappa\left(\frac{C(N)-G(N)}{N}\right)^4,
 \qquad j\ge J.
 \label{eq:dispersion}
\end{equation}
Then
\begin{equation}
 C(n)-\gold n=O\!\left(\frac{n}{(\log n)^{1/4}}\right).
 \label{eq:decay}
\end{equation}
\end{theorem}

\begin{proof}
Let $z=Z(N)$ at a root of order $j$. All labelled values of $X$
lie in $[0,1]$, so its total expected quadratic variation is at most $1$.

Take $H=\lfloor(j-\max(J,\lfloor j/2\rfloor))/8\rfloor$
disjoint four-generation blocks. Each child step lowers order by
at most two, so all nodes in these blocks have order at least
$\max(J,\lfloor j/2\rfloor)$. The difference between
$(C(m)-G(m))/m$ and $Z(m)$ is at most $1/m$.
Write $\varepsilon_j=1/F_{\lfloor j/2\rfloor}$.
At each block start the martingale identity and Jensen's inequality give
\[
 \mathbb E\left[\left(\frac{C(m)-G(m)}m\right)^4\right]
 \ge \bigl((z-\varepsilon_j)_+\bigr)^4.
\]
Summing the assumed four-step variance bounds yields
$1\ge H\kappa((z-\varepsilon_j)_+)^4$.
As $H\asymp j$ and the lower bound gives $z\ge-1/N$,
$|z|=O(j^{-1/4})$. Finally $j\asymp\log N$ proves the result.
\end{proof}

\section{Scalar iteration and inverse rays}
\label{sec:scalar}

The iteration $x\mapsto C(x)$ provides a different view of the Fibonacci
structure. It is distinct from the inner map $T_n$ and from a path through
the two-child recursion. The cap and lower bound imply
\[
 F_{k-1}\le C(x)\le F_k\qquad(F_k<x<F_{k+1}).
\]
Consequently every preimage of a non-Fibonacci value in the open block
$(F_k,F_{k+1})$ lies in the next open block.

\begin{theorem}\label{thm:scalar}
Let $N\ge F_{43}$ and start at a non-Fibonacci integer $x_0>N$.
Iterate $x_{i+1}=C(x_i)$, stopping at the first Fibonacci output or the
first value at most $N$. For any prefix before this stop,
\[
 \sum_i(d(x_i))_+<\frac13,\qquad
 \sum_i|d(x_i)|<\frac13+\frac2N,
\]
where the sums run over the source indices in that prefix.
Every strictly increasing infinite inverse ray
$C(s_{j+1})=s_j$ satisfies
\[
 \sum_j|d(s_j)|<\infty,\qquad s_{j+1}/s_j\longrightarrow\phi.
\]
If the ray is eventually non-Fibonacci, write
$F_{h+j}<s_j<F_{h+j+1}$ after discarding a finite prefix. Then
\[
 \sqrt5(s_j+1)\alpha^{h+j}\ \downarrow\ \lambda\in[1,\phi),
 \qquad s_j/F_{h+j}\longrightarrow\lambda.
\]
\end{theorem}
\begin{proof}
The shifted ratio $\beta(x)=(C(x)+1)/(x+1)$ is strictly greater than
$\alpha$, by the golden floor bound. If $x_0$ is in open block $k$,
all nonterminal steps lower the block order by one. Give the endpoint
$x_r$ the inherited closed block $[F_p,F_{p+1}]$, $p=k-r$.
With $\Lambda_i=\sqrt5(x_i+1)\alpha^{k-i}$ we have exactly
\begin{equation}
 \sum_{i<r}\log\frac{\beta(x_i)}\alpha
 =\log\frac{\Lambda_r}{\Lambda_0}.
 \label{eq:scalar-telescope}
\end{equation}
Binet's formula gives $\Lambda_0>1$. The last source exceeds $N$, so
$x_r>\alpha N-\alpha^2>N/2$. Also
$F_{p+1}\alpha^p\le\phi(1+\alpha^6)/\sqrt5<1$ for $p\ge2$,
whence $\alpha^p<1/x_r<2/N$. A second application of Binet gives
\[
 \Lambda_r\le\phi+\sqrt5\alpha^p+\phi\alpha^{2p+2}
 <\phi+7/N.
\]
Thus the logarithmic sum is less than $\log\phi+7/N$.
The tail bound in Section~\ref{article:maximum-entropy} gives
$C(x)/x<2/3$ at every source. Write $\varepsilon=\beta-\alpha$.
Since $\log(\beta/\alpha)\ge\varepsilon/\beta$ and
$(d(x))_+\le(1+1/N)\varepsilon$,
\[
 \sum_{i<r}(d(x_i))_+
 <(1+1/N)\left(\frac23+\frac1{3N}\right)
                 (\log\phi+7/N)<\frac13.
\]
For the last inequality, $\log\phi<49/100$ and $N\ge F_{43}$ suffice.
The negative part is at most $\alpha^2/x_i$. Reverse geometric summation,
using $x_{i+1}<2x_i/3$, gives $\sum_{i<r}1/x_i<3/N$.
As $3\alpha^2<2$, this proves the absolute bound.

On an inverse ray, once a non-Fibonacci value occurs all later values
are non-Fibonacci, and each edge advances one open block. Reverse a
finite segment ending at a sufficiently large fixed $s_0$ and apply
the preceding bound with $N=s_0$. Monotone partial sums give absolute
summability on the infinite ray. An eventual Fibonacci ray instead
has $s_j=F_{h+j}$, and Binet gives the same conclusion directly.
In either case $d(s_j)\to0$ and the inverse equation gives the ratio limit.

For a non-Fibonacci ray, put
$\Lambda_j=\sqrt5(s_j+1)\alpha^{h+j}$. The inverse equation gives
$\Lambda_{j+1}/\Lambda_j=\alpha/\beta(s_{j+1})<1$, while Binet gives
$\Lambda_j>1$. At every odd block order $p=h+j$, the integer inequality
$s_j+1\le F_{p+1}$ gives
$\Lambda_j\le\phi(1-\alpha^{2p+2})<\phi$. Its limit therefore belongs
to $[1,\phi)$, and the claimed normalization follows from Binet.
\end{proof}

For each fixed $c>0$ and each such ray, absolute summability also gives
\[
 \sum_{\{j:\ d(s_j)\ge c/\sqrt{\log s_j}\}}j^{-1/2}<\infty,
 \qquad
 \#\{j\le J:d(s_j)\ge c/\sqrt{\log s_j}\}=o(\sqrt J).
\]
Indeed, $\log s_j=j\log\phi+O(1)$. It suffices to split the convergent
weighted sum at a fixed index and then let that index tend to infinity.
The theorem concerns every existing inverse ray. Existence from an
arbitrary root and control of the two-child recursion require further
arguments.

\section{Grytczuk's fixed-depth comparison}
\label{sec:grytczuk}

Grytczuk's recurrence \cite{Grytczuk2004}, recorded as \seqnum{A093878}, is
\begin{equation}
 A(1)=A(2)=1,\qquad
 A(n)=A(A(A(n-1)))+A(n-A(A(n-1)))\quad(n\ge3).
 \label{eq:grytczuk}
\end{equation}
It is the $(k,a,b)=(2,0,1)$ case of the family
\[
 A(n)=A(n-a-A^{k}(n-b))+A(A^{k}(n-b))
\]
studied by Isgur and Rahman \cite{IsgurRahman2011}, where $A^k$ denotes
$k$-fold composition.  Their Proposition~1 establishes slow growth and
unboundedness under initial conditions satisfied here, and their
Theorem~2.1 proves generalized Fibonacci anchors.  They credit Grytczuk with the
binary-word interpretation and $A(F_j)=F_{j-1}$ in this case, and explicitly
leave convergence of $A(n)/n$ open for $k\ge2$.  The theorem below settles
the $k=2$ instance of that question and determines the sharp error order.
It makes no assertion for $k>2$.

We derive the required Pascal-word description in Appendix~\ref{app:words}
and use cancellation between a run and its neighboring Pascal word. The
use of Pascal words itself is classical, and
\cite[Sections~1.1, 1.3 and Appendix~B]{Janvresse2005} constructs these blocks
and relates them to the first differences of Conway's sequence.  Here the
Fibonacci interval decomposition and the cancellation between its pieces
give the following quantitative comparison for \eqref{eq:grytczuk}.

\begin{theorem}
The sequence in \eqref{eq:grytczuk} satisfies
\[
 A(n)=\alpha n+O\!\left(\frac n{\sqrt{\log n}}\right).
\]
Moreover,
\[
 \max_{F_{m+2}\le n\le F_{m+3}}|A(n)-\alpha n|
       =\Theta\!\left(\frac{F_{m+1}}{\sqrt m}\right).
\]
In particular, $A(n)/n\to\alpha$, and the global error is not
$o(n/\sqrt{\log n})$.  The lower bound concerns interval maxima.
\end{theorem}

\begin{proof}
For $r\ge1$, define Pascal words by $W(r,0)=0$, $W(r,r)=1$ and
$W(r,j)=W(r-1,j)W(r-1,j-1)$ for $0<j<r$.
For $m\ge2$, the first-difference word on $(F_{m+2},F_{m+3}]$ is
\begin{equation}
 S_m=\prod_{j=\lfloor(m-1)/2\rfloor}^{0}
       1^{\binom{m-2-j}{j}}W(m-1-j,j),
 \label{eq:pascal-segment}
\end{equation}
where the product is in decreasing $j$ order.  The initial endpoint is
$A(F_{m+2})=F_{m+1}$.  Appendix~\ref{app:words} derives these identities.

Put $l=m-1$, $B_j=\binom{l-j}{j}$ and $R_j=\binom{l-j-1}{j}$.
Each run/word pair has length $R_j+B_j$ and $B_j$ ones.  Its excess over
density $\alpha$ is
\[
 e_j=B_j-\alpha(B_j+R_j)=\alpha B_j(\alpha-d_j),
 \qquad d_j=\frac{R_j}{B_j}=\frac{l-2j}{l-j}.
\]
Give $I=l-2j$ probability $B_j/F_{l+1}$.  The polynomial
\[
 P_l(u)=\sum_jB_ju^{l-2j}
 =\frac{r(u)^{l+1}-s(u)^{l+1}}{\sqrt{u^2+4}},\quad
 r(u),s(u)=\frac{u\mathbin{\pm}\sqrt{u^2+4}}2
\]
satisfies $P_l=uP_{l-1}+P_{l-2}$.  Differentiating
$\log P_l(e^t)$ twice at zero gives
$\mathbb EI=l/\sqrt5+O(1)$ and
$\operatorname{Var}I=4l/(5\sqrt5)+O(1)$.
The function $h(x)=2x/(1+x)$ obeys $h(1/\sqrt5)=\alpha$ and
$|h'|\le2$ on $[0,1]$, while $d_j=h(I/l)$.
Cauchy--Schwarz therefore yields
\[
 \sum_j|e_j|=O(F_m/\sqrt m).
\]
This bounds the discrepancy at every complete-pair boundary.

For an unfinished pair, put $p=\alpha^2$ and $t=l-j$.
The identity
\[
 \binom tj p^j(1-p)^{t-j}=\alpha^l B_j
\]
and the standard binomial point-probability bound give
$\max_j B_j=O(F_m/\sqrt m)$, since $t\ge l/2$.
Each pair has length at most $2B_j$, so the same estimate holds for
every prefix.  The starting error satisfies
\[
 F_{m+1}-\alpha F_{m+2}=(-1)^{m+2}\alpha^{m+2}.
\]
Since $n\asymp F_{m+1}$ and $m\asymp\log n$, the upper bound follows.

For the lower bound choose a central index with
$l-2j=l/\sqrt5+O(1)$ on the admissible parity lattice.
Stirling's formula gives $B_j\asymp F_m/\sqrt m$, and
$R_j/B_j\to\alpha$.  The all-one run of length $R_j$ changes
$A(n)-\alpha n$ by $\alpha^2R_j$.  At least one of its two endpoints
has absolute error at least $\alpha^2R_j/2$.  This gives the matching
interval lower bound.
\end{proof}

\section{Discussion and open problems}

The two variable-depth examples have different asymptotic behavior.  The
Campbell orbit is compressed by a ternary scale and is completely classified.
The Cloitre recurrence has a golden Fibonacci skeleton. Exact equality
points, captured collars and the additive profile identity
\eqref{eq:profile} persist across scales, while positive defects create
nontrivial local cycles such as the five-cycle at $196$.

The remaining convergence question is
\[
 \lim_{n\to\infty}\frac{C(n)}n\stackrel{?}=\gold.
\]
The elementary comparison determines the exact equality set and implies
convergence in every sublinear-width Fibonacci neighborhood. The maximum
principle gives the sharper global tail bound, and the convex envelope
records exactly the possible remaining normalized excess.  To obtain the full limit, it suffices to prove sublinear negative
selected jumps when the landing point lies in each fixed interior Fibonacci
subinterval.  The local criterion and the global negative-jump condition
are equivalent for the actual sequence, and localization narrows the
estimate to be proved.  The uniform dispersion condition \eqref{eq:dispersion} is a
separate route to a quantitative rate.
The scalar inverse-ray theorem gives absolute summability along a
different family of paths. Passing from those paths to the selected
two-child recursion is an additional problem, since the scalar estimate
alone does not control its dispersion.

\subsection*{Acknowledgments}
We thank John M. Campbell for the original variable-depth example and the
questions that initiated this investigation.

\subsection*{AI use}
Apart from the variable-depth example of John Campbell, the ideas studied
in this paper are due to the authors, starting with the application of
this idea to Conway's recurrence. The relations with the golden ratio
were first observed in numerical experiments that Beno{\^\i}t Cloitre
designed and ran with Claude Opus 5.5 (Anthropic). Claude was then used to
explore proof strategies and to run further experiments and checks.
ChatGPT (GPT-6 Astra, OpenAI) was used in the same way for the arguments of
Section~\ref{article:maximum-entropy} and Appendix~\ref{app:plateau-collar}.
Haobo Ma and Wenlin Zhang used AI tools for search, double-checking and proofreading.
These were OpenAI Codex (GPT-6 Astra and GPT-6.1 Sol) and ChatGPT
(GPT-6 Pro). All statements, proofs and computations
were checked by the authors, who take full responsibility for the content
and the wording. AI tools also helped with the translation into English.

\subsection*{Code and data availability}
The reproducibility package is available at
\url{https://github.com/the-omega-institute/nested-recurrences}.
It contains independent recurrence evaluators, exact checks of the
Campbell orbit rows, and the finite landing-jump data. The accompanying
source archive contains the programs used for the finite tail certificate,
the small orbit tables, and the profile checks, with reproduction instructions.
The elementary proof of the golden structure in Appendix~\ref{app:golden}
needs only its displayed initial values. The finite computation for the
numerical tail bound covers $[F_{43},F_{45})$. Replaying the small checks
does not require a fresh run of the full prefix. The recorded certificate
and independent evaluator are available as
\href{https://github.com/the-omega-institute/nested-recurrences/tree/974e9d72917b426dfd577503ff705db710b7df34/cloitre-conway/verification/certificate-F45}{versioned verification sources}
in that repository.

\appendix
\section{Deriving the Pascal words from the recurrence}
\label{app:words}

First $A$ is slow-growing and unbounded.  Indeed, if
$s_n=A(A(n-1))$, induction gives $s_n-s_{n-1}\in\{0,1\}$.
Subtracting successive recurrences gives an earlier first difference,
at $n-s_n$ if $s_n=s_{n-1}$, and at $s_n$ otherwise.
The indices are legal since $A(t)\le t-1$ for $t\ge2$, and the base is
$A(3)=2$.  Eventual constancy at $M$ would give
$M=A(A(M))+M$, contradicting positivity.

Let $T(k)$ be the last occurrence of $k$, with $T(0)=0$.  Then
\begin{equation}
 T(1)=2,\qquad T(k)=A(k)+T(k-A(A(k)))\quad(k\ge2).
 \label{eq:last-occurrence}
\end{equation}
The case $k=2$ is direct.  For $k\ge3$, put
$b=A(k)$, $u=A(b)$ and $s=b+T(k-u)$.
Induction and slow growth give $s>T(k-1)$. If $u$ stays fixed, the
earlier endpoint argument advances, whereas an increase in $u$ also
increases $b$ and gives $s=T(k-1)+1$.  At $n=T(k)+1$ the recurrence
gives $A(n-b)=k+1-u$, so $T(k)\ge s$.  If $T(k)>s$, then
$A(s)=A(s+1)=k$, but the recurrence at $s+1$ gives
$u+A(T(k-u)+1)=k+1$, a contradiction.

Write $\mu(k)=T(k)-T(k-1)$.  Subtracting \eqref{eq:last-occurrence},
for $k\ge3$, gives
\[
 \mu(k)=\begin{cases}
 1,&\Delta A(A(k))=1,\\
 \mu(k-A(A(k)))+\Delta A(k),&\Delta A(A(k))=0,
 \end{cases}
\]
where $\Delta A(k)=A(k)-A(k-1)$ and
$\Delta A(A(k))=A(A(k))-A(A(k-1))$.
The initial frequencies are $\mu(1)=2$, $\mu(2)=1$.

Encode each one and its following zeros in $W(r,j)$ by their combined
length, obtaining a word $V(r,j)$ of positive integers.  Thus
$V(1,1)=[1]$, $V(r,1)=[r]$, $V(r,r)=[1]$, and $V$ obeys the same concatenation rule
for $1<j<r$.  Set
\[
 \nu_m=\left(\prod_{j=\lfloor(m-1)/2\rfloor}^{1}
 1^{\binom{m-2-j}{j}}V(m-1-j,j)\right)[2].
\]
Out-of-range words are empty and out-of-range binomial coefficients are zero.
Its gap encoding, replacing $c$ by $10^{c-1}$, is exactly $S_m$.
The word $\nu_m$ has $F_m$ entries and total sum $F_{m+1}$.

For $r\ge2$ and $1\le j\le r$, we use the following replacement identity. In $W(r-1,j-1)$ replace each one by
the integer one, and replace successive zeros by the successive entries
of $V(r-1,j)$ increased by one.  The result is $V(r,j)$.
For $j=1$ and $j=r$ this follows from the edge words.  Otherwise split
the input into $W(r-2,j-1)W(r-2,j-2)$.  Its first part has
$\binom{r-3}{j-1}$ zeros, precisely the number of entries of $V(r-2,j)$.
Splitting the supplied word as $V(r-2,j)V(r-2,j-1)$ and applying induction
proves the assertion.  The row-two cases are checked directly.

We now prove that the frequencies on $(F_{m+1},F_{m+2}]$ are $\nu_m$.
The bases $\nu_2=[2]$, $\nu_3=[1,2]$ follow from
$A(1),\ldots,A(8)=1,1,2,3,3,4,5,5$.
Assume the earlier intervals known.  As $k$ traverses the next interval,
the bits $\Delta A(k)$ traverse $U=S_{m-1}$, and at each one of $U$
the bits $\Delta A(A(k))$ consume the next bit of $Z=S_{m-2}$.
The index $k-A(A(k))$ starts at $F_m$, advances exactly when that
consumed bit is zero, and ends at $F_{m+1}$.  Thus the frequency recurrence
consumes $Q=\nu_{m-1}$ as follows.  A zero of $U$ copies the next $Q$
entry. A one consumes a $Z$ bit and outputs one without consuming $Q$
if that bit is one, and the next $Q$ entry plus one otherwise.

For each decreasing $j\ge1$, put $r=m-1-j$.  The aligned pieces are
\[
\begin{array}{ll}
 U:&W(r-1,j)1^{\binom{r-1}{j-1}},\\
 Z:&1^{\binom{r-2}{j-1}}W(r-1,j-1),\\
 Q:&1^{\binom{r-2}{j}}V(r-1,j).
\end{array}
\]
The $U$ grouping moves the next run after the current Pascal word.
The first $U$ part has $\binom{r-2}{j}$ zeros and
$\binom{r-2}{j-1}$ ones. These consume the initial runs of $Q$ and $Z$
and output $\binom{r-1}{j}$ ones.  The remaining all-one $U$ part
applies the replacement identity to output $V(r,j)$.
Empty pieces handle the top parity cases, and the final zero of $U$
copies the final $[2]$ of $Q$.  This proves $\nu_m$.
Taking sums gives $T(F_{m+2})=F_{m+3}$ and closes the induction.
The gap encoding proves \eqref{eq:pascal-segment} and its endpoint values.

\section{An elementary proof of the golden structure}
\label{app:golden}

We prove the lower bound, equality classification and cap simultaneously.
The proof uses the twenty-four initial values displayed below, and no
large computed prefix. Write $\rho(m)=\{\alpha(m+1)\}$ and
$\Delta(m)=C(m)-G(m)$, reserving $E_k$ for the Fibonacci profiles.
For $a+b=n$,
\begin{equation}
 \begin{aligned}
 \delta(a,b)&=G(a)+G(b)-G(n)
 =-\lfloor\rho(a)+\rho(b)-\alpha\rfloor\in\{-1,0,1\},\\
 \Delta(n)&=\Delta(a)+\Delta(b)+\delta(a,b).
 \end{aligned}
 \label{eq:gold-carry}
\end{equation}
The three carry values $1,0,-1$ correspond to the intervals
$[0,\alpha)$, $[\alpha,1+\alpha)$ and $[1+\alpha,2)$ for
$\rho(a)+\rho(b)$.

\subsection{Initial values and cyclic selection}
\label{app:gold-base}

The following table lists the complete starting orbits through $n=24$.
A semicolon separates the transient from the cycle, and ranks start at zero.
Every transition $u\mapsto v$ is the substitution $C(u)=n-v$ using an
already known value. The selected rank is
$\tau+((h-\tau)\bmod p)$, where $h=C(n-1)$ and $p$ is the cycle length.
\begin{center}\small
\begin{longtable}{rrlrr}
\caption{The initial orbits for $C$.}\label{tab:gold-orbits}\\
\toprule $n$&$h$&Transient; (cycle)&$g(n)$&$C(n)$\\\midrule
\endfirsthead
\toprule $n$&$h$&Transient; (cycle)&$g(n)$&$C(n)$\\\midrule
\endhead
3&1&$\,;\ (2)$&2&2\\
4&2&$\,;\ (3,2)$&3&3\\
5&3&$\,;\ (4,2)$&2&3\\
6&3&$5;\ (3,4)$&3&4\\
7&4&$6,3,5;\ (4)$&4&5\\
8&5&$7,3,6,4;\ (5)$&5&5\\
9&5&$8,4;\ (6,5)$&5&6\\
10&6&$9,4;\ (7,5)$&7&7\\
11&7&$10,4,8;\ (6,7)$&6&7\\
12&7&$11,5,9,6,8;\ (7)$&7&8\\
13&8&$12,5,10,6,9,7;\ (8)$&8&8\\
14&8&$13,6,10,7;\ (9,8)$&9&9\\
15&9&$14,6,11;\ (8,10)$&8&10\\
16&10&$15,6,12,8,11;\ (9,10)$&10&11\\
17&11&$16,6,13,9,11;\ (10)$&10&12\\
18&12&$17,6,14,9,12,10;\ (11)$&11&12\\
19&12&$18,7,14,10;\ (12,11)$&12&13\\
20&13&$19,7,15,10,13;\ (12)$&12&13\\
21&13&$20,8,16,10,14,12;\ (13)$&13&13\\
22&13&$21,9,16,11,15,12;\ (14,13)$&13&14\\
23&14&$22,9,17,11,16,12;\ (15,13)$&15&15\\
24&15&$23,9,18,12;\ (16,13)$&13&15\\\bottomrule
\end{longtable}\end{center}
In particular $h\ge\tau$ in every row. The comparisons are
\begin{center}\small
\begin{tabular}{r|rrrrrrrrrrrr}
\toprule
$n$&1&2&3&4&5&6&7&8&9&10&11&12\\\midrule
$C(n)$&1&1&2&3&3&4&5&5&6&7&7&8\\
$G(n)$&1&1&2&3&3&4&4&5&6&6&7&8\\
$U(n)$&1&1&2&3&3&4&5&5&6&7&8&8\\\midrule
$n$&13&14&15&16&17&18&19&20&21&22&23&24\\\midrule
$C(n)$&8&9&10&11&12&12&13&13&13&14&15&15\\
$G(n)$&8&9&9&10&11&11&12&12&13&14&14&15\\
$U(n)$&8&9&10&11&12&13&13&13&13&14&15&16\\\bottomrule
\end{tabular}\end{center}
Here $U(1)=1$. The bounds $21/34<\alpha<34/55$ determine all
floors in this table, since the two rational products have the same floor
for every multiplier from $2$ to $25$. The equality indices are exactly
$Z\cap[1,24]=\{1,2,3,4,5,6,8,9,11,12,13,14,21,22,24\}$.

\begin{proof}[Proof of Lemma~\ref{lem:foundations}]
Assume the earlier terms are positive and $C(m)\le m$, with the strict
bound $C(m)\le m-1$ for $m\ge2$. Then $T_n$ maps $[1,n-1]$ into
$[2,n-1]$, since its image at $1$ is $n-1$ and its value at $x\ge2$ is
at least $n-x+1$. Thus the children are earlier indices, at least one exceeds
$1$, and their values sum to at most $n-1$. Adding the inductive
half-bounds gives $C(n)\ge\lceil n/2\rceil$.

The values through $5$ initialize $C(n)\le3n/4$. For $n\ge6$,
the earlier upper bound gives $x_i\ge(n+3)/4>2$ for $i\ge1$.
Since $h\ge2$, both $a=x_h$ and $b=C(x_{h-1})$ are at least $2$,
so the upper bounds add. The values through $9$ initialize
$C(n)\ge3n/5$. For $n\ge10$ we have $a,x_{h-1}\ge4$ and
$C(y)\ge3$ for $y\ge4$, so both children are at least $3$ and the
lower bounds add as well.

For $n\ge6$ these bounds successively give
\[
 x_i\ge\frac{n+3}{4}\ (i\ge1),\quad
 x_i\le\frac{17n-9}{20}\ (i\ge2),\quad
 x_i\ge\frac{29n+27}{80}\ (i\ge3).
\]
The integer interval containing the tail from $x_3$ has at most
$(39n+17)/80$ elements. Distinctness before the first repetition gives
\[
 \tau\le2+\frac{39n+17}{80}
 =\frac{39n+177}{80}\le\frac{3(n-1)}5\le h\qquad(n\ge25).
\]
The middle inequality is $9n\ge225$. Table~\ref{tab:gold-orbits}
handles $3\le n\le24$, proving cyclic selection for every $n\ge3$.
\end{proof}

\subsection{Rotation arithmetic}
The Fibonacci recurrence gives
\begin{equation}
 \alpha F_j-F_{j-1}=(-1)^{j-1}\alpha^j.
 \label{eq:fib-rotation}
\end{equation}
Hence, for $j\ge3$,
\[
 \rho(F_j)=\alpha+(-1)^{j-1}\alpha^j,\quad
 \rho(F_j+1)=2\alpha-1+(-1)^{j-1}\alpha^j,
 \quad \rho(F_j-1)=\alpha^j\quad(j\text{ odd}).
\]
At $11,24,25,59$ the fractions are respectively
$12\alpha-7$, $25\alpha-15$, $26\alpha-16$ and $60\alpha-37$.
Each lies in $(0,\alpha)$, as follows from $37/60<\alpha<5/8$.
Consequently, if $z\in Z$, $z\ge8$, then $\rho(z)>\alpha$ exactly
when $z=F_k$ with odd $k\ge7$. Also
$\rho(z)\le\alpha+\alpha^7$, so two such zeros have
$\rho(a)+\rho(b)\le2\alpha+2\alpha^7<1+\alpha$.
Thus their carry is nonnegative.

\begin{lemma}\label{lem:one-sided-rotation}
If $k\ge5$ is odd and $0<t<F_{k+2}$ is an integer, then
$\{\alpha t\}\ge1-\alpha^k$ if and only if $t=F_{k+1}$.
\end{lemma}
\begin{proof}
First, for $j\ge3$, $0<t<F_j$ and any integer $s$,
\begin{equation}
 |\alpha t-s|\ge\alpha^{j-1},
 \label{eq:fib-best}
\end{equation}
with equality only at $(t,s)=(F_{j-1},F_{j-2})$.
Cassini's identity makes $(F_j,F_{j-1})$ and
$(F_{j-1},F_{j-2})$ an integer basis. Write $(t,s)$ as their
linear combination with coefficients $u,v$. The error is
$\alpha^{j-1}|u\alpha-v|$. If $u\ge1$, the bound $t<F_j$
forces $v\le-1$, and if $u\le-1$, positivity forces $v\ge1$.
Both give a factor greater than $1$. If $u=0$, $v\ge1$ and equality
requires $v=1$.

For $t<F_{k+1}$, apply \eqref{eq:fib-best} with $s=\lceil\alpha t\rceil$.
The proposed fraction would force $t=F_k$, whose error has the wrong
sign. At $F_{k+1}$ the fraction is $1-\alpha^{k+1}$ and qualifies.
Finally let $t=F_{k+1}+w$, $0<w<F_k$. The same bound gives
$\alpha^{k-1}\le\{\alpha w\}\le1-\alpha^{k-1}$.
Subtracting $\alpha^{k+1}$ does not cross zero and gives a fraction
strictly below $1-\alpha^k$.
\end{proof}

\subsection{Capture under prefix hypotheses}
Assume for now, for $x<n$, the three assertions
\[
 C(x)\ge G(x),\qquad C(x)=G(x)\Longrightarrow x\in Z,
 \qquad C(x)\le U(x).
 \tag{$\mathcal P$}
\]
The function $U$ is nondecreasing with increments $0$ or $1$.
For an anchor $A=F_k\ge5$, $B=F_{k-1}$, the following cones hold
at arguments $x<n$.
\begin{equation}
 \begin{array}{c|c}
 x<A&B+x-A+1\le C(x)\le B\\
 x\ge A&B\le C(x)\le B+x-A.
 \end{array}
 \label{eq:strict-cones}
\end{equation}
The upper bounds follow from $U$ and the lower bound on the second
row follows from $G$. For the strict first bound write $x=A-u$.
If $u\ge2$,
$\alpha(A-u+1)-(B-u+1)=\alpha^2(u-1)+(-1)^{k-1}\alpha^k>0$.
For $u=1$ and odd $k$, $G(A-1)=B$. For even $k\ge6$,
$G(A-1)=B-1$, but $A-1\notin Z$ (including the four exceptions),
so the zero implication in $(\mathcal P)$ supplies the additional unit.

Let $I$ be the closed interval with endpoints $A,n-B$, and write
its endpoints as $\ell,h$. The cones make $I$ invariant and give
\[
 T_n(\ell-u)\in[\ell,h+u-1],\qquad
 T_n(h+u)\in[\ell-u,h]\qquad(u\ge1)
\]
whenever the arguments belong to the domain. Thus the integer distance
to $I$ decreases in at most two steps until it is zero. This proof does
not evaluate $C(A)$, so the anchor may equal $n$.

Write $q=F_j\le n<F_{j+1}=q+p$, $p=F_{j-1}$, $r=F_{j-2}$.
For $n\ge8$, applying capture at $q$ and $p$ gives, for every cycle point,
\begin{equation}
 \max(n-p,p)=:\ell\le x\le h:=\min(q,n-r).
 \label{eq:cycle-rectangle}
\end{equation}
Since selection is cyclic, the oriented children obey
\begin{equation}
 p\le a=g(n)\le q,\qquad r\le b=n-a\le p.
 \label{eq:child-rectangle}
\end{equation}
A cycle predecessor $y$ of $a$ satisfies $C(y)=b$. Moreover, every cycle
point is at most $n-G(\ell)$.
Put $s=F_{j-3}$, $t=F_{j-4}$. The two child caps are
\begin{equation}
 C(a)\le\min(a-s,p),\qquad C(b)\le\min(b-t,r).
 \label{eq:child-caps}
\end{equation}
Adding separately the affine and constant bounds proves
$C(n)\le\min(n-r,q)=U(n)$.

We strengthen the induction at the changes of slope
$K_j=F_j+F_{j-2}$ by requiring
\begin{equation}
 C(K_j)<U(K_j)\ (j\ge6),\qquad
 C(K_j\pm1)<U(K_j\pm1)\ (j\ge8).
 \label{eq:strict-knees}
\end{equation}
These statements follow before the lower bound at the current index.
Indeed, let $A=K_{j-1}=p+s$, $B=K_{j-2}=r+t$ and $n=K_j+v$.
Equality throughout \eqref{eq:child-caps} would force
\[
\begin{array}{c|c|c}
v&U(n)&(a,b)\\\hline
0&q&(A,B)\\
1&q&(A+1,B)\text{ or }(A,B+1)\\
-1&q-1&(A-1,B)\text{ or }(A,B-1).
\end{array}
\]
For $v=0,1$ both constant caps must be attained, and for $v=-1$ both
affine caps must be attained. In each row both children would saturate their
own caps. The seed is $C(K_6)=C(11)=7<8$. For $j\ge7$ the central
row contradicts the earlier strict center at $K_{j-1}$, and for $j\ge8$
each neighbor row forces a strict center of order $j-1$ or $j-2$ to
saturate. This proves the simultaneous propagation of
\eqref{eq:strict-knees}.

\subsection{The simultaneous induction}
The initial table proves $(\mathcal P)$ through $24$ and initializes
\eqref{eq:strict-knees}. At $n\ge25$ assume both statements at all
earlier indices. The preceding arguments establish the cap and any
strict-knee assertion at $n$. We prove $\Delta(n)\ge0$ and its zero
implication. In its Fibonacci block, $j\ge8$, $p\ge13$, $r\ge8$.

If both child defects are positive, \eqref{eq:gold-carry} gives
$\Delta(n)\ge1$. If exactly one is zero, the parent defect is nonnegative
and can vanish only with carry $-1$. Write $z$ for the zero argument
and $w$ for the other. Then $\rho(z)>\alpha$, so $z=F_k$ for odd
$k\ge7$, and $\{\alpha(w+1)\}\ge1-\alpha^k$.
If $z=a$, then $b<a$, and Lemma~\ref{lem:one-sided-rotation} excludes
$b+1\le F_k$. If $z=b=p=F_k$, the lemma forces
$a+1=q=F_{k+1}$, whence $n=F_{k+2}-1\in Z$.
If $z=b=r=F_k$, then $a+1\in[F_{k+1}+1,F_{k+2}+1]$.
The lemma excludes the interior below $F_{k+2}$, and the last two
fractions are $\alpha^{k+2}$ and $\alpha+\alpha^{k+2}$, both below
$1-\alpha^k$ since $\alpha^k+\alpha^{k+2}<\alpha^2$.

Suppose both child defects vanish. Their carry is nonnegative by the
rotation arithmetic. Apart from $11,24,25,59$, the possibilities are
\[
 a\in\{p,p+1,q\}\cup\{q-1:j\text{ odd}\},\quad
 b\in\{r,r+1,p\}\cup\{p-1:j-1\text{ odd}\}.
\]
We enumerate the zero-carry pairs.

For consecutive anchors write $a=F_{k+1}+s$, $b=F_k+t$, with
$s,t\in\{-1,0,1\}$ and offset $-1$ permitted only at odd orders.
The offsets cannot both be $-1$. If $s+t=0$ or $1$, $n\in Z$.
If $s+t=2$, the fractional-part sum is
$4\alpha-2+(-1)^{k-1}\alpha^{k+2}<\alpha$, giving carry $1$.
If $s+t=-1$ and $k$ is odd, again $n\in Z$. If $k$ is even, then
$(s,t)=(-1,0)$ and the sum $\alpha-\alpha^{k+2}$ gives carry $1$.

For anchors two ranks apart write $a=q+s$, $b=r+t$, where
$s\in\{-1,0\}$ and $t\in\{0,1\}$. A predecessor $y$ satisfies
$y\ge n-p=2r+s+t$ but $C(y)=r+t$. Since $\alpha>3/5$,
\[
 \frac35(2r+s+t+1)-(r+t+1)=\frac{r+3s-2t-2}{5}
 \ge\frac{r-7}{5}>0,
\]
contradicting $C(y)\ge G(y)$.

For the same anchor put $p=F_k$, $k=j-1$, and $B=2p-r=K_k$.
The pair $(p+1,p-1)$, allowed only for odd $k$, has carry $1$ because
$2\alpha-1+2\alpha^k<\alpha$. The other pairs require
\[
\begin{array}{c|c|c}
(a,b)&C(y)&y\\\hline
(p,p)&p&B\\
(p+1,p)&p&B\text{ or }B+1\\
(p,p-1)&p-1&B-1.
\end{array}
\]
The upper bounds on $y$ come from capture, and the lower bounds follow
from $U(y)\le y-F_{k-2}$. Each value saturates $U(y)$, contradicting
\eqref{eq:strict-knees} for $k\ge8$. At the remaining order $j=8$,
$(13,12)$ has sum $25\in Z$, $(13,13)$ contradicts the strict center
$K_7=18$, and $(14,12)$ has carry $1$. The pair $(14,13)$ at $n=27$
is excluded by the earlier values giving
$14\mapsto18\mapsto15\mapsto17\mapsto15$, in which $14$ is transient.

The exceptional children occur only in the following products of lists.
\[
\begin{array}{c|c|c}
q&a&b\\\hline
21&13,14,21&11\\
34&24,25&13,14,21\\
55&34,35,55&24,25\\
89&59&34,35,55\\
144&89,90,144&59.
\end{array}
\]
At $q\ge233$, $b\ge89$, so this list is exhaustive. Remove $(13,11)$,
already in the base. Of the remaining twenty pairs,
\[
 (25,14),\ (35,25),\ (59,35),\ (90,59)
\]
have carry $1$. The pairs $(14,11),(34,25),(35,24)$ have sums
$25,59,59$ in $Z$.
The other thirteen are excluded as follows. All these floors are fixed
by $55/89<\alpha<89/144$.

For $(21,11),(55,24),(55,25),(144,59)$ the triples
$(n,\ell,G(\ell))$ are
\[
 (32,19,12),\quad(79,45,28),\quad(80,46,29),\quad(203,114,71).
\]
A predecessor would have $C(y)\ge G(\ell)>b$.
For three further pairs, two successive cycle predecessors give
\[
\begin{array}{c|c|c|c|c|c}
n&(a,b)&[\ell,h]&y&C(z)&z\\\hline
37&(24,13)&[21,24]&21&16&24\\
38&(25,13)&[21,25]&21&17&25\\
93&(59,34)&[55,59]&55&38&59.
\end{array}
\]
Here $G(22)=14$ or $G(56)=35$ forces $y$, and then the affine cap
$U(z)\le z-8$ or $z-21$ forces $z=a$. Its required value is
$G(a)+1$, contradicting the zero-child assumption. This does not assume
$C(59)$ in advance. The pair $(24,14)$ at $n=38$ is excluded by
$24\mapsto23\mapsto23$, using $C(24)=C(23)=15$.
For the next three pairs, the refined upper capture bound gives
\[
\begin{array}{c|c|c|c|c}
n&(a,b)&\ell&n-G(\ell)&\text{possible }y\\\hline
45&(24,21)&24&30&29,30\\
46&(25,21)&25&30&29,30\\
114&(59,55)&59&77&76,77.
\end{array}
\]
The lower bounds on $y$ follow from $U(y)\ge b$. But $29=K_8$ and
$76=K_{10}$, so the strict center and neighbor inequalities exclude
these values.

The last two pairs require the prescribed depth, not just cyclicity.
If $T_n(H)\subseteq J$ and $T_n(J)\subseteq H$, an orbit starting
in $H$ is in $H$ at even ranks. The needed intervals are
\[
\begin{array}{c|c|c|c|c|c}
n&H&J&G(\min H)&U(\max J)&h\\\hline
58&[35,57]&[1,36]&22&23&36\\
148&[90,147]&[1,92]&56&58&92.
\end{array}
\]
The displayed comparisons prove the two inclusions. Since
$57\notin Z$, $G(57)=35$ and $U(57)=36$, the depth at $58$ is $36$.
Similarly, $147\notin Z$, $G(147)=91$ and $U(147)=92$ force depth $92$
at $148$. Thus $a\ge35$ or $90$, excluding $(34,24)$ and $(89,59)$.
This completes all cases, hence the simultaneous induction.

\subsection{The reverse inclusion and landing}
We now have $(\mathcal P)$ at every index. The equalities
$G(F_j)=U(F_j)=F_{j-1}$ and
$G(F_j+1)=U(F_j+1)=F_{j-1}+1$ give the anchors.
The strict cone and $U(F_j-1)=F_{j-1}$ give the left neighbors for
$j\ge5$. When $j$ is odd, these equal $G(F_j-1)$.
The smaller cases are in the initial table.

The table gives $C(11)=G(11)$ and $C(24)=G(24)$. At $n=25$ the orbit
$24\mapsto10\mapsto18\mapsto13\mapsto17\mapsto13$
has transient $3$, period $2$ and depth $15$, selecting $13$.
Thus $C(25)=C(13)+C(12)=16=G(25)$.
At $59$, the bounds and $58\notin Z$ force depth $C(58)=37$.
The intervals $H=[37,58]$ and $J=[1,35]$ alternate, since $C(37)\ge24$
by its strict lower bound, $C(x)\ge G(38)=24$ for $x\ge38$, and
$C(x)\le U(35)=22$ for $x\le35$. Odd depth selects $a\le35$.
Capture also gives $34\le a\le38$, so $a=34$ or $35$.
Both give $C(59)=37=G(59)$ using $C(25)=16$ and $C(24)=15$.
This proves the exact equality set.

At $n=F_k$, $k\ge6$, capture with $A=F_{k-1}$ makes the interval
with endpoints $A,n-F_{k-2}$ a singleton. Every orbit therefore
reaches $A$, which is fixed by the anchor identity. Finally the lower
bound gives $\liminf C(n)/n\ge\alpha$, and the Fibonacci subsequence
gives the reverse inequality.

\section{Proofs of the left plateau and of the right collar}
\label{app:plateau-collar}

Both proofs use Theorems~\ref{thm:golden} and~\ref{thm:cap}, the
anchors~\eqref{eq:fib-values}, and the capture and child localization
in Theorem~\ref{thm:capture} and~\eqref{eq:child-rectangle}, with inherited
labels at Fibonacci endpoints. By Lemma~\ref{lem:foundations}, $g(n)$
lies on the cycle $\Gamma_n$ reached from $n-1$.
The function $G$ is nondecreasing and $1$-Lipschitz, and $U$ is
nondecreasing.

\subsection*{The left plateau}

\begin{proof}[Proof of Theorem~\ref{article:left-plateau}]
Put $L_k(u)=F_{k-1}-C(F_k-u)$ for $0\le u\le F_{k-2}$. We prove the
stronger invariant
\begin{equation}
\label{app:left-invariant}
\begin{aligned}
 L_k(u)&=0 &&(0\le u\le w_k),\\
 L_k(w_k+1)&=1,\\
 1\le L_k(u)&\le u-w_k-1 &&(w_k+2\le u\le F_{k-2}).
\end{aligned}
\end{equation}
A direct evaluation of the recursion through $55$ gives the following
vectors $(L_k(u))_{0\le u\le F_{k-2}}$.
\[
\begin{array}{c|l}
 k& (L_k(u))_{0\le u\le F_{k-2}}\\ \hline
 8&(0,0,0,1,1,2,3,4,5)\\
 9&(0,0,0,0,1,1,2,2,3,5,6,6,7,8)\\
 10&(0,0,0,0,1,1,1,1,2,2,3,4,6,5,6,8,9,9,10,11,12,13)
\end{array}
\]
These $45$ values verify \eqref{app:left-invariant} for $k=8,9,10$.
The levels $9$ and $10$ initialize the induction below, whose parents
satisfy $n\ge F_{10}=55>47$.

Fix $k\ge11$ and suppose that \eqref{app:left-invariant} holds at the
levels $k-1$ and $k-2$. Write
\[
 A=F_{k-1},\qquad B=F_{k-2},\qquad P=F_{k-3},\qquad Q=F_{k-4},
\]
so that $A=B+P$ and $B=P+Q$. Let $n=F_k-u$ with $0\le u\le B$.
With the inherited label $k-1$ for the closed block $[A,F_k]$, the
selected children satisfy $g(n)\in[B,A]$ and $n-g(n)\in[P,B]$.
In particular $g(n)\in[A-u,A]$.

The interval $[A-u,A]$ is invariant under $T_n$. Indeed, if $y=A-v$
with $0\le v\le u$, then $C(y)\le U(A)=B$, and
$C(y)\ge G(y)\ge G(A)-v=B-v\ge B-u$. Hence $A-u\le T_n(y)\le A$.
Since $g(n)$ lies on $\Gamma_n$, the whole cycle lies in $[A-u,A]$.
Capture at level $k-1$ also gives $y\ge B$ on $\Gamma_n$, and for a
cycle predecessor $z$ of $y$ the bounds $G(B)=P$ and $U(A)=B$ give
$n-y=C(z)\in[P,B]$. Every cycle point can therefore be written
$y=A-v$ with
\begin{equation}
\label{app:left-cycle-domain}
 0\le v\le u,\qquad v\le P,\qquad 0\le u-v\le Q.
\end{equation}
In these coordinates the cycle map and the additivity at the selected
point read
\begin{equation}
\label{app:left-cycle-map}
 v\longmapsto u-L_{k-1}(v),
 \qquad
 L_k(u)=L_{k-1}(v)+L_{k-2}(u-v).
\end{equation}

Set $w=w_{k-1}$ and $z=w_{k-2}$, so that $z\ge2$. If $0\le u\le w$,
every accessible value $L_{k-1}(v)$ is zero, every cycle point is sent
to $u$, and the cycle is the fixed point $v=u$. Then
\eqref{app:left-cycle-map} gives $L_k(u)=0$.

Suppose that $u\ge w+2$. No cycle point has $v\le w$, since it would be
sent to $u$. If $u>P$, this image is excluded by
\eqref{app:left-cycle-domain}. If $u\le P$, the invariant at level
$k-1$ gives $1\le L_{k-1}(v)\le u-w-1$ for $w<v\le u$, so after
reaching $u$ the orbit stays in $[w+1,u-1]$ and cannot return to a
point with $v\le w$. A cycle point $v=u$ would have a predecessor with
zero defect, hence with offset at most $w$, which is excluded.
Therefore $w+1\le v\le u-1$ on $\Gamma_n$. For $u=w+2$ the cycle is the
fixed point $v=w+1$, whose complementary offset is $1\le z$, and
\eqref{app:left-cycle-map} gives
\begin{equation}
\label{app:left-next}
 L_k(w+2)=1.
\end{equation}
For $u\ge w+3$, put $b=u-v\ge1$. If $v=w+1$, then $b\ge2$, and the
invariant at level $k-2$ with $z\ge2$ gives $L_{k-2}(b)\le b-2$, so
$L_k(u)\le1+(b-2)=u-w-2$. If $v\ge w+2$, the bounds
$L_{k-1}(v)\le v-w-1$ and $L_{k-2}(b)\le b-1$ give the same estimate.
In both cases $L_{k-1}(v)\ge1$, hence
\begin{equation}
\label{app:left-far}
 1\le L_k(u)\le u-w-2\qquad(w+3\le u\le B).
\end{equation}

It remains to determine $L_k(w+1)$. The defect at level $k-1$ vanishes
through $w$ and equals $1$ at $w+1$, so the only cycle of
\eqref{app:left-cycle-map} in the allowed range is $\{w,w+1\}$. In the
original coordinates its vertices are $\ell=A-w-1$ and $r=A-w=\ell+1$.
The invariant at level $k-1$ gives the barriers
\begin{equation}
\label{app:left-barriers}
 y<r\ \Longrightarrow\ C(y)\le B-1,
 \qquad
 y\ge r\ \Longrightarrow\ C(y)\ge B.
\end{equation}
For $B\le y<r$ the first one is the positivity of the defect at level
$k-1$, and for $y<B$ it follows from $C(y)\le U(B)=P<B$. For
$r\le y\le A$ the second one is the plateau at level $k-1$, and for
$y>A$ it follows from $C(y)\ge G(A)=B$. Since $n=A+B-w-1$, the
barriers give
\[
 y\ge r\ \Longrightarrow\ T_n(y)\le\ell,
 \qquad
 y\le\ell\ \Longrightarrow\ T_n(y)\ge r.
\]
The orbit starts at $n-1=r+B-2\ge r$ and alternates between the two
sides. Once it has reached the cycle, it visits $r$ exactly at even
ranks. By \eqref{app:left-next}, the depth at this parent is
$h_n=C(n-1)=A-L_k(w+2)=A-1$. This is not circular, since
\eqref{app:left-next} was obtained from the two previous profiles
without the value at $w+1$. If $A$ is odd, the depth is even, the
selected point is $r$ with offset $w$, and $L_k(w+1)=0$. If $A$ is
even, the selected point is $\ell$ with offset $w+1$, and
$L_k(w+1)=1$.

Consequently $w_k=w_{k-1}+\mathbf 1_{\{F_{k-1}\text{ odd}\}}$. Since
$F_j$ is even exactly when $3$ divides $j$, the formula
$w_k=\lfloor(2k-9)/3\rfloor$ satisfies this recurrence and the initial
values. Equations \eqref{app:left-next} and \eqref{app:left-far} give
the value $1$ just beyond the new plateau and the affine bound of
\eqref{app:left-invariant}, whether the plateau advances or not. This
completes the induction.

Finally, if $m<F_{k-1}$, then $C(m)\le U(F_{k-1})=F_{k-2}<F_{k-1}$.
The invariant excludes every index of $[F_{k-1},F_k-w_k)$, and the
plateau starts at $F_k-w_k$. This gives the first occurrence.
\end{proof}

\subsection*{The right collar}

Let $n=F_k+t\ge47$ with $0\le t\le F_{k-1}$, and write
$g(n)=F_{k-1}+s$ with inherited labels. Capture and child localization
give $0\le s\le t$, $s\le F_{k-2}$ and $t-s\le F_{k-3}$, and
Theorem~\ref{thm:capture} gives $E_k(t)=E_{k-1}(s)+E_{k-2}(t-s)$, where all defects
are nonnegative by the cap. Before entering the captured interval, the
orbit may leave the range where $P_{k-1}$ is used. For that purpose
put $\widetilde P_{k-1}(u)=C(F_{k-1}+u)-F_{k-2}$ whenever
$1\le F_{k-1}+u\le n-1$. In the coordinate $y=F_{k-1}+u$ the map $T_n$
becomes $u\mapsto t-\widetilde P_{k-1}(u)$, and the monotonicity of
$U$ and $G$ together with the anchors gives the barriers
\begin{equation}
\label{app:collar-barriers}
 \widetilde P_{k-1}(u)\le0\quad(u<0),
 \qquad
 \widetilde P_{k-1}(u)\ge1\quad(u>0).
\end{equation}

\begin{lemma}[Transmission of a first defect]
\label{app:first-defect}
Let $1\le t\le F_{k-3}$ with $F_k+t\ge47$, and suppose that
$E_j(u)=0$ for $0\le u<t$ and $j\in\{k,k-1,k-2\}$. If $E_k(t)>0$, then
\[
 E_{k-1}(t)=0,\qquad E_k(t)=E_{k-2}(t)>0,\qquad g(F_k+t)=F_{k-1},
 \qquad F_{k-1}\equiv t\pmod2.
\]
\end{lemma}

\begin{proof}
Write $A=F_{k-1}$ and $B=F_{k-2}$. The hypotheses force $k\ge10$, so
$B>1$. Every cycle lies in the coordinate interval $[0,t]$, and the
selected point is cyclic.

Suppose first that $E_{k-1}(t)>0$. Since $P_{k-1}(t)\ge1$, the integer
$E_{k-1}(t)$ lies in $[1,t-1]$. For $0\le u<t$ the map is
$u\mapsto t-u$, and the image of $t$ is $E_{k-1}(t)$. Thus $0$ has no
preimage in $[0,t]$, and the only preimage of $t$ there is $0$. Neither
endpoint is cyclic, so $0<s<t$, and both terms of the profile identity
vanish, contrary to $E_k(t)>0$. Hence $E_{k-1}(t)=0$.

The map on $[0,t]$ is now the reflection $u\mapsto t-u$. If $s>0$,
then $E_{k-1}(s)=0$ and $E_{k-2}(t-s)=0$, which again contradicts
$E_k(t)>0$. Hence $s=0$ and $E_k(t)=E_{k-2}(t)>0$.

It remains to find the parity of the selected rank. The initial
coordinate is $u_0=n-1-A=B+t-1>t$. By \eqref{app:collar-barriers},
before its first entry into $[0,t]$ the orbit lies above $t$ at even
ranks and below $0$ at odd ranks. An entry at a point of $(0,t)$ stays
in $(0,t)$ under the reflection and cannot reach the endpoint cycle
$\{0,t\}$. An entry into that cycle therefore happens at $0$ at an odd
rank or at $t$ at an even rank, and in both cases every visit to $0$
has odd rank. Since the selected offset is $0$, the depth $h_{F_k+t}$
is odd. The collar at $t-1$ gives $h_{F_k+t}=C(F_k+t-1)=A+t-1$, hence
$A\equiv t\pmod2$.
\end{proof}

\begin{proof}[Proof of Proposition~\ref{article:log-collar}]
For $t\ge0$ put $K_t=8+6t$. We prove by induction on $t$ that
$E_k(u)=0$ for $0\le u\le t$ and every $k\ge K_t$. The case $t=0$ is
the anchor $C(F_k)=F_{k-1}$. Let $t\ge1$, assume the assertion for
$t-1$, fix $k\ge K_t$ and suppose that $E_k(t)>0$. We apply
Lemma~\ref{app:first-defect} at the levels $k$, $k-2$ and $k-4$. The
lowest profile index involved is $k-6\ge K_{t-1}$, which gives the
collar hypotheses. The offset hypotheses follow from
$t\le F_{6t+1}\le F_{k-7}$, and the lowest parent level is
$k-4\ge10$, so every parent exceeds $47$. The lemma gives
$E_{k-2}(t)>0$, then $E_{k-4}(t)>0$, and
\[
 F_{k-1}\equiv F_{k-3}\equiv F_{k-5}\equiv t\pmod2.
\]
The indices $k-1$, $k-3$ and $k-5$ represent the three residue classes
modulo $3$, so exactly one of these Fibonacci numbers is even, a
contradiction. Hence $E_k(t)=0$, and the smaller offsets are covered by
the induction hypothesis.
\end{proof}

\section{The Campbell formula and its orbit rows}\label{app:campbell}
Define $B(1)=1$ and define $B(n)$ by the right side of
\eqref{eq:campbell-formula} for $n\ge2$. Its scale is unique because the respective
half-open intervals partition the eligible even and odd indices. Directly,
\begin{equation}\label{camp:eq:bounds}
1\le B(n)\le n-1,\qquad B(n)\equiv n+1\pmod2\quad(n\ge2).
\end{equation}
The following four rules hold on closed intervals, and the formulas agree
at every shared endpoint. In each row $t$ is a power of $3$.
\begin{center}
\begin{tabular}{lll}
\toprule
Parity of $x$ & Interval & $B(x)$\\
\midrule
even & $2t\le x\le4t$ & $x-t$\\
even & $4t\le x\le6t$ & $3t$\\
odd & $3t\le x\le5t$ & $2t$\\
odd & $5t\le x\le9t$ & $x-3t$\\
\bottomrule
\end{tabular}
\end{center}
At $x=6t$, for instance, the next even scale $3t$ gives the same value
$3t$, and at $x=9t$ the next odd scale gives $6t$.

Fix $n$ and put $U_n(x)=n-B(x)$ on $1\le x\le n-1$. By
\eqref{camp:eq:bounds}, this interval is invariant. Write $x_0=n-1$ and
$x_{j+1}=U_n(x_j)$.

\begin{lemma}\label{camp:lem:orbit}
Let $s\ge3$ be a power of $3$ in the interval associated with $n$.
The orbit is given by the tables below. In particular, $U_n^6(n-1)=U_n^4(n-1)$.
If $n$ is even then $x_4=B(n)$, and if $n$ is odd then $x_5=B(n)$.
\end{lemma}

In the first table set
\[
A=\min(n-2s/3,5s/3+1),\qquad B=\min(n-2s/3,2s).
\]
Rows are closed for convenience, and are intersected with $2s\le n<6s$
or $3s\le n<9s$, respectively. Adjacent rows agree at shared endpoints.

\begin{center}
\textbf{Even $n$}\par\smallskip
\setlength{\tabcolsep}{3pt}
\resizebox{\textwidth}{!}{%
\begin{tabular}{c c c c c c c c}
\toprule
& Interval for $n$ & $x_1$ & $x_2$ & $x_3$ & $x_4$ & $x_5$ & $x_6$\\
\midrule
E1 & $[2s,3s+1]$ & $s+1$ & $n-2s/3-1$ & $A$ & $n-s$ & $B$ & $n-s$\\
E2 & $[3s+1,10s/3]$ & $n-2s$ & $7s/3$ & $n-4s/3$ & $n-s$ & $2s$ & $n-s$\\
E3 & $[10s/3,4s]$ & $n-2s$ & $n-s$ & $2s$ & $n-s$ & $2s$ & $n-s$\\
E4 & $[4s,5s+1]$ & $n-2s$ & $3s$ & $n-2s$ & $3s$ & $n-2s$ & $3s$\\
E5 & $[5s+1,6s]$ & $3s+1$ & $n-2s-1$ & $n-2s$ & $3s$ & $n-2s$ & $3s$\\
\bottomrule
\end{tabular}}
\end{center}

\begin{center}
\textbf{Odd $n$}\par\smallskip
\setlength{\tabcolsep}{3pt}
\resizebox{\textwidth}{!}{%
\begin{tabular}{c c c c c c c c}
\toprule
& Interval for $n$ & $x_1$ & $x_2$ & $x_3$ & $x_4$ & $x_5$ & $x_6$\\
\midrule
O1 & $[3s,4s+1]$ & $s+1$ & $n-2s/3-1$ & $5s/3+1$ & $n-s$ & $2s$ & $n-s$\\
O2 & $[4s+1,13s/3]$ & $n-3s$ & $10s/3$ & $n-7s/3$ & $n-s$ & $2s$ & $n-s$\\
O3 & $[13s/3,5s]$ & $n-3s$ & $n-s$ & $2s$ & $n-s$ & $2s$ & $n-s$\\
O4 & $[5s,6s+1]$ & $n-3s$ & $4s$ & $n-3s$ & $4s$ & $n-3s$ & $4s$\\
O5 & $[6s+1,7s]$ & $3s+1$ & $n-2s-1$ & $n-3s$ & $4s$ & $n-3s$ & $4s$\\
O6 & $[7s,8s+1]$ & $3s+1$ & $n-2s-1$ & $n-3s$ & $n-3s$ & $n-3s$ & $n-3s$\\
O7 & $[8s+1,9s]$ & $3s+1$ & $n-2s-1$ & $5s+1$ & $n-3s$ & $n-3s$ & $n-3s$\\
\bottomrule
\end{tabular}}
\end{center}

\begin{proof}
Use the four branch rules with $t=s/3,s,3s$, all powers of $3$. The
parities alternate for even $n$, starting with odd $x_0$, and for odd $n$
all indices $x_j$ are even. All listed arguments remain in $[1,n-1]$.

We detail E1, which involves the two minima. Here
$5s/3\le2s-1\le n-1\le3s$, so the odd linear branch at scale $s/3$
gives $B(n-1)=n-s-1$ and $x_1=s+1$. Since
$2s/3\le s+1\le4s/3$, we get $B(x_1)=2s/3+1$ and
$x_2=n-2s/3-1$. This odd index lies in $[s,3s]$, and hence
\[
B(x_2)=\max(2s/3,n-5s/3-1),\qquad x_3=A.
\]
Both arguments defining $A$ are at least $4s/3$, and
$A\le5s/3+1\le2s$. The even plateau at scale $s/3$ therefore gives
$B(x_3)=s$ and $x_4=n-s$. This odd index lies in $[s,3s]$, so
$B(x_4)=\max(2s/3,n-2s)$ and $x_5=B$. Finally,
$4s/3\le B\le2s$ gives $B(x_5)=s$ and $x_6=n-s$.

The remaining rows follow by the same substitutions. To illustrate the
change of scale in O7, $x_1=3s+1$ gives $B(x_1)=2s+1$, so
$x_2=n-2s-1\in[6s,7s-1]$. The even linear branch at scale $3s$
gives $B(x_2)=n-5s-1$ and $x_3=5s+1$. The even plateau at scale $s$
then gives $x_4=n-3s$, which lies in $[5s+1,6s]$ and is fixed by $U_n$.
Each of the other rows uses the displayed branch intervals directly.
\end{proof}

\subsection{Induction for the formula}
We prove $b(n)=B(n)$ by strong induction. Direct evaluation through $n=8$
gives $(1,1,2,3,2,3,4,5)$, as required. For $n\ge9$, assume equality
at every smaller index. Then $T_n=U_n$ throughout $[1,n-1]$. Since this
interval is invariant, their entire starting orbits agree.

Let $d=b(n-1)=B(n-1)$. Equation~\eqref{camp:eq:bounds} gives
$d\equiv n\pmod2$. For even $n\ge10$, the table gives $d\ge4$, since
$d=n-s-1\ge n/2-1\ge4$ in E1, $d=2s\ge6$ in E2--E4, and
$d=n-3s-1\ge2s\ge6$ in E5. For odd $n$, it gives $d\ge5$, since
$d=n-s-1\ge2s-1\ge5$ in O1, $d=3s\ge9$ in O2--O4, and
$d=n-3s-1\ge3s\ge9$ in O5--O7.

Since $x_6=x_4$, determinism implies $x_{j+2}=x_j$ for every $j\ge4$.
For even $n$, $d$ is even and at least four, so
$b(n)=x_d=x_4=B(n)$. For odd $n$, $d$ is odd and at least five, so
$b(n)=x_d=x_5=B(n)$. This closes the induction.

The same lemma proves $T_n^6(n-1)=T_n^4(n-1)$ for $n\ge9$. For
$n=2,3,4,5,6,7,8$, direct calculation gives $x_4=x_6$ equal to
$1,2,3,4,3,4,5$, respectively. The bound four is attained at $n=24$, with the orbit
\[
23\longmapsto10\longmapsto17\longmapsto16\longmapsto15
\longmapsto18\longmapsto15\longmapsto18\longmapsto\cdots.
\]
The first six displayed values are distinct. This finishes the proof of
Theorem~\ref{thm:campbell}.

\subsection{Scaling and periods}
\begin{corollary}
For every $n\ge2$, $b(3n)=3b(n)$. Moreover,
\[
\liminf_{n\to\infty}\frac{b(n)}n=\frac25,\qquad
\limsup_{n\to\infty}\frac{b(n)}n=\frac34.
\]
\end{corollary}
\begin{proof}
Multiplying $n$ by $3$ preserves parity, sends its scale $s$ to $3s$,
and multiplies either expression in~\eqref{eq:campbell-formula} by $3$.
On even indices, $b(n)/n\in[1/2,3/4]$, and on odd indices it lies in
$[2/5,2/3]$. The values $n=5\cdot3^k$ attain $2/5$, and
$n=4\cdot3^k$ attain $3/4$, proving the two limits.
\end{proof}

For even $n\ge4$, the eventual period is exactly two, since the two
phases have opposite parity. At $n=2$ the orbit is fixed. For odd $n$,
the tables give
\[
x_4=\begin{cases}
n-s,&3s\le n\le5s,\\
4s,&5s\le n\le7s,\\
n-3s,&7s\le n<9s,
\end{cases}
\qquad x_5=\max(2s,n-3s).
\]
Consequently the prescribed orbit is eventually fixed exactly when
$n=3s$ or $7s\le n<9s$. At the other odd indices its period is two.
The small scale $s=1$ agrees by direct calculation at $n=3,5,7$.

\bigskip\hrule\bigskip
\noindent 2020 \emph{Mathematics Subject Classification}:
11B37; 11B39; 68W30.

\noindent \emph{Keywords:} nested recurrences, meta-recursions, Fibonacci
numbers, Hofstadter sequences, Pascal words.

\end{document}